\documentclass[12pt]{amsart}
\usepackage{amscd}
\usepackage{amssymb,amsmath,amsthm}
\usepackage[all]{xy}

\theoremstyle{plain}
\newtheorem{thm}{Theorem}[section]
\newtheorem{lem}[thm]{Lemma}

\newtheorem{cor}[thm]{Corollary}
\newtheorem{defn-lem}[thm]{Definition-Lemma}

\newtheorem{prop}[thm]{Proposition}

\theoremstyle{definition}
\newtheorem{defn}[thm]{Definition}
\newtheorem{ex}[thm]{Example}
\newtheorem{rem}[thm]{Remark}

\def\md #1#2#3#4#5 {\left(
                        \begin{matrix}
             #1 & #2 \\
             #3 & #4
                        \end{matrix}
                      \right)- #5}

\def\ma #1#2#3#4 {\left(
                        \begin{matrix}
             #1 & #2 \\
             #3 & #4
                        \end{matrix}
                      \right)}

\def\Index{\operatorname{Index}}

\def \Ker {\operatorname{Ker}}

\def\tsr{\operatorname{tsr}}
\def\Ad{\operatorname{Ad}}
\def\id{\operatorname{id}}
\def\Aut{\operatorname{Aut}}

\newcommand{\mc}{\mathcal}
\newcommand\bigzero{\makebox(0,0){\text{\huge0}}}

\numberwithin{equation}{section}

\begin{document}
\title[Rokhlin property for inclusions]
       {Rokhlin property and approximate representability for inclusions of
       index-finite type}

\begin{abstract}
Let $P\subset A$ be an inclusion of unital $C\sp*$-algebras of index-finite
type with a fixed conditional expectation $E:A\to P$.  We introduce the
weak tracial Rokhlin property and its dual notion, weak tracial approximate
representability, using positive contractions in central sequence algebras; in
particular, the definitions apply to projectionless algebras.  Under natural
simplicity, finiteness, and outerness hypotheses, we prove that these properties
are exchanged by the Jones--Watatani basic construction and its dual conditional
expectation.  We show that the associated Rokhlin contractions produce
tracially large injective completely positive order-zero maps and that our
inclusion-theoretic definition recovers the weak tracial Rokhlin property for
finite-group actions.  The duality also yields inclusions of index-finite type with the
weak tracial Rokhlin property which are not isomorphic to fixed-point
inclusions arising from actions of ordinary finite groups.  Explicit actions on
an infinite-type UHF algebra, a simple monotracial AF algebra, and $\mc Z$
separate approximate, tracial approximate, and weak tracial approximate
representability. Finally, by showing that the inclusion
$P\hookrightarrow A$ is tracially sequentially-split by order zero, we obtain
permanence of tracial $\mc{Z}$-absorption, $\mc{Z}$-stability, strict
comparison, tracial $m$-comparison, tracial $m$-almost divisibility, and
tracial nuclear dimension from $A$ to $P$.
\end{abstract}

\author { Hyun Ho \, Lee and Hiroyuki Osaka}

\address {Department of Mathematics\\
          University of Ulsan\\
         Ulsan, 44610 South Korea }
\email{hadamard@ulsan.ac.kr}

\address{Department of Mathematical Sciences\\
Ritsumeikan University\\
Kusatsu, Shiga 525-8577, Japan}
\email{osaka@se.ritsumei.ac.jp}


\date{}
\thanks{The first author's research was partially supported by Basic Science Research Program through the National Research Foundation of Korea(NRF) funded by the Ministry of Education(2018R1D1A1B07050924). 
The second author's research was partially supported by KAKENHI Grant Number JP17K05285 and JP20k03644.}
\maketitle

\section{Introduction}

Rokhlin-type freeness is one of the principal mechanisms by which a symmetry
of a $C\sp*$-algebra can preserve regularity.  For an outer action
$\alpha:G\curvearrowright A$ of a finite group, this symmetry is encoded not
only by the action but also by the inclusions of index-finite type
\[
        A^{\alpha}\subset A\subset A\rtimes_{\alpha}G
\]
and their canonical conditional expectations.  This observation raises two
questions. Can one formulate an intrinsic Rokhlin theory for inclusions of
index-finite type that is compatible with Jones--Watatani duality? How far does
the resulting class extend beyond ordinary finite-group actions?

The strict Rokhlin property for inclusions of index-finite type was introduced
by Osaka, Kodaka, and Teruya in \cite{OKT:Rokhlin}; Osaka and Teruya developed
its tracial version in \cite{OT1}, with a correction in \cite{OT2}.  Both
definitions, like Izumi's strict and Phillips' tracial Rokhlin properties for
finite-group actions \cite{Izumi:finite,Izumi:Rokhlin,Phillips:tracial}, use
projections and hence do not reach stably finite projectionless algebras such
as the Jiang--Su algebra $\mc Z$.  In the group-action setting, positive
contractions in central sequence algebras, with Cuntz comparison controlling
the complementary defect, provide the appropriate replacement for Rokhlin
projections \cite{GHS}.  The inclusions associated with finite-group actions
therefore point to a corresponding positive-contraction theory for general
inclusions.

This passage is not merely formal.  A Rokhlin contraction must retain an exact
relation with the Jones projection, and its Cuntz-small defect must be
transported among the central sequence algebras of the smaller algebra, the
larger algebra, and the basic construction.  At the level of maps, this
requires injective completely positive order-zero maps in place of corner
homomorphisms.  These requirements shape our intrinsic definition.

Let $P\subset A$ be an inclusion of unital $C\sp*$-algebras of index-finite
type, and fix a conditional expectation $E:A\to P$ of index-finite type.  We
introduce the
\emph{weak tracial Rokhlin property} for $E$ and its dual notion, the
\emph{weak tracial approximate representability}.  Both notions are formulated
in ultrapowers and require the complementary defect to be Cuntz subequivalent
to an arbitrarily prescribed nonzero positive element.  The compatibility
with the basic construction is encoded by
\[
              (\Index E)e^{1/2}e_Pe^{1/2}=e,
\]
where $e$ is the Rokhlin contraction.  Thus the definition is intrinsic to
the conditional expectation and remains meaningful without projections.

The definition first passes a basic consistency test.  If
$\alpha:G\curvearrowright A$ is an action of a finite group on a simple,
separable, unital, infinite-dimensional $C\sp*$-algebra, then $\alpha$ has the
weak tracial Rokhlin property if and only if the averaging expectation
$E:A\to A^{\alpha}$ is of index-finite type and has the weak tracial Rokhlin property
(Theorem~\ref{T:action-inclusion}).  Hence the inclusion-theoretic definition
recovers the dynamical one.  More generally, a Rokhlin contraction gives a
tracially large injective c.p.c. order-zero map
\[
                  A\longrightarrow P_{\omega}
\]
and, when $E$ is outer, the existence of these maps characterizes the
property (Theorem~\ref{T:orderzerofromA}).  This order-zero formulation is
also the mechanism that later yields permanence results.

Our principal results are the two duality theorems.  Let
$B=C\sp*\langle A,e_P\rangle$ be the basic construction and let
$\widehat{E}:B\to A$ be the dual conditional expectation.  Under the natural
simplicity and outerness assumptions, $E$ has the weak tracial Rokhlin property
if and only if $\widehat{E}$ is weakly tracially approximately representable
(Theorem~\ref{T:tracialduality1}).  If $A$ is finite, the complementary
equivalence also holds: $E$ is weakly tracially approximately representable if
and only if $\widehat{E}$ has the weak tracial Rokhlin property
(Theorem~\ref{T:tracialduality2}).  Consequently, the two notions form a dual
pair under the Jones--Watatani basic construction, extending the familiar
Rokhlin/approximate-representability duality beyond projection-based settings.
The proofs require transferring Cuntz comparison through both expectations
and the basic construction while retaining the exact Jones--Watatani
relations.  This is closely related to the realization of crossed products by
saturated finite-dimensional Hopf $C\sp*$-algebra actions as basic constructions
due to Szyma\'nski and Peligrad \cite{SP:saturated}.

The duality is not confined to inclusions of ordinary group type.  Starting
with an outer weakly tracially approximately representable action of a finite
non-abelian group $G$, we prove that the inclusion
\[
                  A\subset A\rtimes_{\alpha}G,
\]
with its canonical expectation
\[
                  F:A\rtimes_{\alpha}G\longrightarrow A
\]
has the weak tracial Rokhlin property.  Nevertheless, this inclusion, with its
specified expectation, is not isomorphic to a fixed-point inclusion arising
from an action of an ordinary finite group
(Proposition~\ref{P:nonabelian-dual-example}).  The obstruction is the
dual Hopf $C\sp*$-algebra $C(G)$.  Every group Hopf $C\sp*$-algebra
$C\sp*(H)$ is cocommutative, whereas the coproduct on $C(G)$ is
cocommutative exactly when $G$ is abelian.  Thus, for non-abelian $G$,
$C(G)$ cannot be isomorphic to any $C\sp*(H)$ as a Hopf $C\sp*$-algebra.
The inclusion framework therefore detects finite quantum symmetry that is
invisible at the level of ordinary finite-group actions.

This should be compared with \cite{LOT:Rokhlin-index}, where the present
authors and Teruya construct non-group-type Rokhlin and tracial Rokhlin
inclusions using quantum symmetries on the purely infinite strongly
self-absorbing algebras $\mathcal O_2$ and $\mathcal O_\infty$.  In those
examples, non-integer index separates the inclusions from classical
fixed-point models and, in particular cases, even from fixed-point models for
finite-dimensional Hopf $C\sp*$-algebras.  The constructions in the present
paper instead start from stably finite algebras: the infinite-type UHF algebra
$M_{|G|^\infty}$, a simple monotracial AF algebra, and the projectionless
algebra $\mc Z$.  The UHF algebra and $\mc Z$ are strongly self-absorbing
\cite{TW:ssa}.  The resulting inclusions have integer index $|G|$, and their
non-group-type obstruction comes from the dual Hopf algebra $C(G)$ rather
than from nonintegrality of the index.

Three explicit constructions realize genuinely different levels of
representability.  For every finite non-abelian group $G$, the product-type
action
\[
 \alpha_g=\bigotimes_{n=1}^{\infty}\Ad(\lambda_g)
       \quad\text{on}\quad M_{|G|^\infty},
\]
where $\lambda$ is the left regular representation, is approximately
representable and pointwise outer.  Both assertions are proved directly in
Example~\ref{E:UHF-AR}.  Finite tensor products of the $\lambda_g$ give exact
implementing representations in the sequence algebra, while the uniform
spectral separation of $\lambda_g$ from the scalar unitaries prevents the
phase-adjusted partial implementers from being norm-Cauchy when $g\ne1_G$.

At the next level, we construct an action on a simple monotracial AF algebra
which is tracially approximately representable but not approximately
representable.  At the weakest level, equivariant prime dimension-drop maps
produce a pointwise outer action of an arbitrary finite non-abelian group on
$\mc Z$ which is weakly tracially approximately representable but is neither
tracially approximately representable nor approximately representable.
Thus the three examples exhibit the hierarchy from approximate to tracial
approximate to weak tracial approximate representability, with both reverse
implications failing.  Applying the duality theorem gives three corresponding
families of non-group-type inclusions with the weak tracial Rokhlin property.
In particular, the $\mc Z$ example proves that the weak positive-contraction
theory is strictly broader than both projection-based notions, even in the
monotracial projectionless setting.

As consequences of the order-zero formulation, we obtain a unified package of
regularity permanence results.  We first prove that $P$ is simple whenever
$A$ is simple and $E$ has the weak tracial Rokhlin property.  We then show
that the canonical inclusion $P\hookrightarrow A$ is tracially sequentially
split by order zero in the sense of \cite{LeeOsaka4}.  The permanence theorems
of \cite{LeeOsaka4} and \cite{LeePermanenceII} consequently pass tracial
$\mc Z$-absorption, $\mc Z$-stability, strict comparison, tracial
$m$-comparison, tracial $m$-almost divisibility, and tracial nuclear dimension
from $A$ to $P$.  By duality, analogous conclusions hold for weak tracial
approximate representability and the basic construction.  These applications
demonstrate the utility of the theory; the duality and the strict-separation
examples are the main new contributions.

The paper is organized as follows.  Section~\ref{S:notation} reviews
inclusions of index-finite type, basic constructions, and order-zero maps.
Section~\ref{S:Rokhlin} introduces the weak tracial Rokhlin property, its
order-zero characterization, and its relation with finite-group actions.
Section~\ref{S:approxrepresentability} develops weak tracial approximate
representability.  Section~\ref{S:tracial} proves the two duality theorems and
then constructs the UHF, AF, and $\mc Z$ examples and the resulting
inclusions beyond ordinary group actions.  The final section develops the tracially
sequentially-split order-zero formulation and records the resulting permanence
theorems.

\section{Preliminaries and notation}\label{S:notation}
We recall inclusions of index-finite type and their conditional expectations,
the Jones--Watatani basic construction, and completely positive order-zero
maps. We also fix the notation used throughout the paper.
\begin{defn}[Watatani]
Let $P\subset A$ be an inclusion of unital $C\sp*$-algebras and let
$E:A\to P$ be a conditional expectation. A finite family
$\{(u_j,v_j)\}_{j=1}^n\subset A\times A$ is a \emph{quasi-basis} for $E$ if
\[
 x=\sum_{j=1}^n u_jE(v_jx)=\sum_{j=1}^nE(xu_j)v_j
 \quad\text{for every }x\in A.
\]
If such a family exists, then $E$ is said to be of \emph{index-finite type},
and its Watatani index is
\[
                    \Index E=\sum_{j=1}^n u_jv_j.
\]
\end{defn}
We say that the inclusion $P\subset A$ is \emph{of index-finite type} if it
admits a faithful conditional expectation of index-finite type.  Whenever a
particular expectation $E:A\to P$ is part of the data, it is fixed and assumed
to be faithful and of index-finite type.  We reserve the term \emph{Watatani index} for the
element $\Index E$ and do not use it as a synonym for this class of
inclusions.
Whenever a quasi-basis exists, one can choose a self-adjoint quasi-basis
$\{(u_j,u_j^*)\}_{j=1}^n$. Moreover, $\Index E$ is a nonzero positive element
of $Z(A)$. In particular, if $A$ is simple, then $\Index E$ is a positive
scalar.

We next recall the basic construction. Since all conditional expectations in
this paper are of index-finite type, the reduced and maximal basic
constructions agree.
\begin{defn}[Watatani]\cite[Chapter II]{W:index}
Assume that $E:A\to P$ is faithful. Let $\mc{E}_E$ be the Hilbert $P$-module
completion of $A$ for the $P$-valued inner product
$\langle x,y\rangle_P=E(x^*y)$. Let $\mc{L}(\mc{E}_E)$ denote the algebra of
adjointable operators on $\mc{E}_E$, let
$\lambda:A\to\mc{L}(\mc{E}_E)$ be the representation by left multiplication,
and let $\eta_E:A\to\mc{E}_E$ be the canonical map. The \emph{Jones
projection} $e_P$ is defined by
\[ e_P(\eta_E(x))=\eta_E(E(x)).\]
The \emph{basic construction} is the $C\sp*$-algebra
\[C\sp* \langle A, e_P \rangle =\{\sum_{i=1}^n \lambda(x_i)e_P \lambda(y_i)\mid x_i, y_i \in A, n\in \mathbb{N} \}. \]
When $E$ is of index-finite type, there is a dual conditional expectation
$\widehat{E}:C\sp*\langle A,e_P\rangle\to A$ satisfying, for $x,y\in A$,
\[
 \widehat{E}(\lambda(x)e_P\lambda(y))=(\Index E)^{-1}xy.
\]
Moreover, $\widehat{E}$ is also of index-finite type and faithful, and
$\Index\widehat{E}=\Index E$.
Furthermore, by \cite[Lemma 2.11]{W:index},
\begin{equation}
e_P x e_P= E(x)e_P
\end{equation}
for any $x\in A$. 
\end{defn}
Unless otherwise stated, all conditional expectations are faithful. We
suppress $\lambda$ and write $xe_Py$ for
$\lambda(x)e_P\lambda(y)$.

We will also use faithfulness after passing to norm ultrapowers.  If
$F:C\to D$ is of index-finite type, the Pimsner--Popa inequality gives a
constant $c>0$, independent of $x\in C$, such that
\[
                         F(x^*x)\ge c x^*x.
\]
Applying this estimate coordinatewise shows that
$F_{\omega}(X^*X)\ge cX^*X$ for $X\in C_{\omega}$.  In particular,
$F_{\omega}$ is faithful.  This applies both to $E_{\omega}$ and to the
ultrapower $\widehat E_{\omega}$ of the dual conditional expectation.

We will use the following structure theorem.
\begin{thm}\cite[Theorem 3.3]{Izumi:inclusion}\label{T:structure}
Let $P\subset A$ be an inclusion of unital $C\sp*$-algebras of index-finite type. Then 
\begin{enumerate}
\item If $A$ is simple, then $P$ is a finite direct sum of simple closed two-sided ideals.
\item If $P$ is simple, then $A$ is a finite direct sum of simple closed two-sided ideals. 
\end{enumerate}
Moreover, if $A=\bigoplus_i A_i$, then $A_i=Az_i$ for a central projection
$z_i\in Z(A)$.
\end{thm}

\begin{defn}\cite[Definition 1.3]{WZ:orderzero}
Positive elements $a,b$ in a $C\sp*$-algebra are \emph{orthogonal}, written
$a\perp b$, if $ab=0$. A completely positive map $\phi:A\to B$ has
\emph{order zero} if
\[ a\perp b \Longrightarrow \phi(a)\perp \phi(b).\] 
\end{defn}
We refer to such a map simply as an order-zero map. We will repeatedly use the
following structure theorem from \cite{WZ:orderzero}.
\begin{thm}[Winter and Zacharias]\label{T:orderzero}
Let $A$ and $B$ be $C\sp*$-algebras, and let $\phi : A \to B$ be an order-zero map. Set $C :=C\sp*(\phi(A)) \subset B$.

Then, there are a positive element  $h \in \mathcal{M}(C)\cap C'$ with $\| h\|= \| \phi \|$ and a $*$-homomorphism 
\[\pi_{\phi}: A \to \mathcal{M}(C)\cap \{h\}'\]
such that 
\[\pi_{\phi}(a)h=\phi(a) \,\text{for $a\in A$}.\]
If $A$ is unital, then $h=\phi(1_A)\in C$. 
\end{thm}
The following simple observation will be used repeatedly. 
\begin{lem}\label{L:simpleorderzero}
Let $A$ be a simple unital $C\sp*$-algebra and $\phi:A \to B$ be a nonzero order-zero map. Then $\phi$ is injective. 
\end{lem}
\begin{proof}
Write $\phi(\cdot)=h\pi_{\phi}(\cdot)$ as in Theorem~\ref{T:orderzero}. Consider $a\in A$ and $b\in \Ker \phi$. Then 
\[ \begin{split}
\phi(ab)&=h\pi_{\phi}(ab)\\
&=h\pi_{\phi}(a)\pi_{\phi}(b)\\
&=\pi_{\phi}(a)h\pi_{\phi}(b)\\
&=\pi_{\phi}(a)\phi(b)=0
\end{split}
\]
Therefore $ab,ba\in\Ker\phi$, and $\Ker\phi$ is a closed two-sided ideal.
Since $\phi$ is nonzero and $A$ is simple, $\Ker\phi=0$.
\end{proof}
\begin{lem}\label{L:lemma}
Let $P\subset A$ be an inclusion of unital $C\sp*$-algebras of index-finite
type, and let
\[
          \psi:C\sp*\langle A,e_P\rangle\longrightarrow D
\]
be an order-zero map into a $C\sp*$-algebra $D$, with $\psi(e_P)=e$.
Write $\psi(\cdot)=h\pi_{\psi}(\cdot)$ as in
Theorem~\ref{T:orderzero}. If $a\in A$ and $\psi(ae_P)e=0$, then
$\pi_{\psi}(E(a)e_P)e=0$.
\end{lem}
\begin{proof}
Since $\psi(1)=h$, $e=\psi(e_P)\le h$. Note that  $\psi(ae_P)e=0$ implies that $\pi_{\psi}(ae_P)he=0$. 
Then \[
\begin{split}
0 &\le \pi_{\psi}(ae_P)ee\pi_{\psi}(ae_P)^*\\
&\le \pi_{\psi}(ae_P)e^{1/2}ee^{1/2}\pi_{\psi}(ae_P)^*\\
&\le  \pi_{\psi}(ae_P)e^{1/2}he^{1/2}\pi_{\psi}(ae_P)^*\\
&\le \pi_{\psi}(ae_P)he\pi_{\psi}(ae_P)^*=0 \quad (\text{$h$ commutes with $\psi(e_P)$.})
\end{split}
\]
Thus \[\pi_{\psi}(ae_P)e=0.\]
So \[\pi_{\psi}(e_P)\pi_{\psi}(ae_P)e=\pi_{\psi}(e_Pae_P)e = \pi_{\psi}(E(a)e_P)e=0. \]
\end{proof}
\section{The weak tracial Rokhlin property for an inclusion of unital $C\sp*$-algebras of index-finite type}\label{S:Rokhlin}
Fix a free ultrafilter $\omega$ on $\mathbb N$. For a $C\sp*$-algebra $A$,
let $\ell^{\infty}(\mathbb N,A)$ be the algebra of bounded sequences in $A$
and let
\[
c_{\omega}(\mathbb N,A)
=\{(a_n)_n:\lim_{n\to\omega}\|a_n\|=0\}.
\]
The norm ultrapower is
\[
A_{\omega}=\ell^{\infty}(\mathbb N,A)/c_{\omega}(\mathbb N,A),
\]
and $\pi_{\omega}$ denotes the quotient map. We write
$a=[(a_n)_n]$ for an element of $A_{\omega}$. We identify $A$ with the
constant sequences and write its central sequence algebra as
\[A_{\omega} \cap A'.\]
After changing a representing sequence outside a set belonging to $\omega$,
every nonzero projection $p\in A_{\omega}$ may be represented by nonzero
projections $p_n\in A$. Likewise, a nonzero positive element
$a\in A_{\omega}$ may be represented by nonzero positive elements $a_n\in A$
whose norms are uniformly bounded away from zero.

We begin with an observation based on \cite{OKT:Rokhlin,OT1}.
\begin{prop}\label{P:equivalentconditions}
Let $P\subset A$ be an inclusion of unital $C\sp*$-algebras, where $A$ is
simple and separable, and let $E:A\to P$ be a conditional expectation of
index-finite type. For a projection $e\in A_{\omega}\cap A'$, set
$g=(\Index E)E_{\omega}(e)$. Then the following are equivalent:

\begin{enumerate}
\item $(\Index E)ee_Pe=e$;
\item $g$ is a projection in $P_{\omega}$;
\item $ge=e=eg$.
\end{enumerate}
\end{prop}
\begin{proof}
$(1)\Rightarrow(2)$. It suffices to prove that $g^2=g$. Indeed,
 \[\begin{aligned}
g^2e_P&= (\Index E) E_{\omega}(e)  (\Index E) E_{\omega}(e)  e_P\\
&= (\Index E)^2 E_{\omega}(e) e_P e e_P\\
&=(\Index E)^2e_Pe e_P e e_P\\
&=(\Index E) e_Pe e_P\\
&=ge_P
\end{aligned}\] 
The map $P_{\omega}\ni z\mapsto ze_P$ is injective, since
\[
             \widehat E_{\omega}(ze_P)=(\Index E)^{-1}z.
\]
Therefore $g^2=g$.

$(2)\Rightarrow(3)$. Since $A$ is simple, the Pimsner--Popa inequality gives
\[ E(x^*x) \ge \frac{x^*x}{(\Index E)^2}.\]
It follows that 
 \[\begin{aligned}
E_{\omega}(e)&\ge \frac{e}{(\Index E)^2},\\
g=(\Index E)E_{\omega}(e) &\ge \frac{e}{(\Index E)}.
\end{aligned}\] 
Therefore \[(1-g)\frac{e}{\Index E}(1-g)=0.\] 
\[(e(1-g))^*(e(1-g))=0.\] Thus $e(1-g)=0$. 

$(3)\Rightarrow(1)$. Let $f=(\Index E)ee_Pe$. Then
\[\begin{aligned}
f^2&= (\Index E) ee_Pe (\Index E) e e_Pe\\
&= (\Index E)^2 e e_P e e_Pe\\
&=(\Index E) e (\Index E)E_{\omega}(e)e_P e\\
&=(\Index E) e g e_P e\\
&=(\Index E)ee_Pe=f
\end{aligned}\] 
Thus $f$ is a projection with $f\le e$. Since
$\widehat{E}_{\omega}(e-f)=0$, the faithfulness of
$\widehat{E}_{\omega}$ implies that $e=f$.
\end{proof}

For positive elements $a,b\in A$, we write $a\lesssim b$ if there is a
sequence $(x_n)_n$ in $A$ such that $x_nbx_n^*\to a$. If $p$ is a projection,
then $p\lesssim a$ if and only if the hereditary subalgebra generated by $a$
contains a projection Murray--von Neumann equivalent to $p$. We refer to
\cite{Cu,Ro:UHF1,Ro:UHF2} for further details.

Proposition~\ref{P:equivalentconditions} isolates the two structural
components of the tracial Rokhlin property.  Suppose that $A$ is simple and
has property~(SP).  In the projection-valued formulation of \cite{OT1}, for
every nonzero positive element $z\in A_{\omega}$ one finds a projection
$e\in A_{\omega}\cap A'$ such that, with
\[
                     g=(\Index E)E_{\omega}(e),
\]
the element $g$ is a projection and $1-g\lesssim z$.  By
Proposition~\ref{P:equivalentconditions}, the first requirement is equivalent
to the Jones--Watatani compatibility relation
\[
                         (\Index E)ee_Pe=e,
\]
and it implies $ge=e$.  Thus the definition combines an exact compatibility
condition with the basic construction and a tracial largeness condition,
expressed by the Cuntz-smallness of $1-g$.

The Jones--Watatani relation has a natural positive-contraction analogue,
\[
                   (\Index E)e^{1/2}e_Pe^{1/2}=e.
\]
This formulation does not require projections and therefore remains
meaningful for projectionless algebras.  It is also compatible with the
dynamical picture: for the averaging expectation associated to a finite-group
action,
\[
             (\Index E)E_{\omega}(e)
             =\sum_{g\in G}\alpha_{g,\omega}(e),
\]
so the second condition below is precisely the usual requirement that the
Rokhlin contractions exhaust the unit up to a Cuntz-small defect.  These
considerations lead to the following definition.

\begin{defn}\label{D:Rokhlin}
Let $P\subset A$ be an inclusion of unital $C\sp*$-algebras, where $A$ is
simple, and let $E:A\to P$ be a conditional expectation of index-finite type.
We say that $E$ has the \emph{weak tracial Rokhlin property} if, for every
nonzero positive element $a\in A_{\omega}$, there is a nonzero positive
contraction $e\in A_{\omega}\cap A'$ such that
\begin{enumerate}
\item $(\Index E) e^{1/2}e_Pe^{1/2}=e$, 
\item $1-(\Index E)E_{\omega}(e) \lesssim a$.
\end{enumerate} 
\end{defn}
\begin{rem}
If $A$ is not assumed to be simple, one may add the following nondegeneracy
condition:
 \begin{itemize}
 \item[(3)] the map $ A \ni x \mapsto xe= e^{1/2}xe^{1/2} \in A_{\omega}$ is injective. 
 \end{itemize}
If $A$ is simple, this condition follows from
Lemma~\ref{L:simpleorderzero}.
 \end{rem}

Sometimes, it is convenient for us to consider the following equivalent condition for (1) in Definition~\ref{D:Rokhlin}.

\begin{prop}\label{P:equivalentcondition}
Let $P\subset A$ be an inclusion of unital $C\sp*$-algebras, where $A$ is
simple, and let $E:A\to P$ be a conditional expectation of index-finite type.
For a positive contraction $e\in A_{\omega}\cap A'$, the following are
equivalent:
\begin{enumerate}
\item $(\Index E) e^{1/2}e_P e^{1/2}=e$.
\item $(\Index E)E_{\omega}(e)e^{1/2}=e^{3/2}$.
\end{enumerate}
\end{prop}
\begin{proof}
Suppose that $(\Index E)e^{1/2}e_Pe^{1/2}=e$. Then, for $a\in A$,
\[\begin{aligned}
ae^{3/2}&= (\Index E) \widehat{E}_{\omega}(e_Pae^{3/2})\\
&= (\Index E)  \widehat{E}_{\omega}(e_Pa (\Index E)ee_Pe^{1/2})\\
&=(\Index E)^2 \widehat{E}_{\omega}(E_{\omega}(ae)e_P e^{1/2}) \\
&=(\Index E)E_{\omega}(ae)e^{1/2}\\
\end{aligned}\] 
Taking $a=1_A$ gives
$e^{3/2}=(\Index E)E_{\omega}(e)e^{1/2}$.

Conversely, let $f=(\Index E)e^{1/2}e_Pe^{1/2}$. Then
\[\begin{aligned}
f^2&= (\Index E)^2 e^{1/2}e_P e e_P e^{1/2}\\
&= (\Index E)^2 e^{1/2} E_{\omega}(e)e_P e^{1/2}\\
&=(\Index E) (\Index E)E_{\omega}(e)e^{1/2}e_P e^{1/2}\\
&=(\Index E) e^{3/2} e_P e^{1/2}\\
&=ef.
\end{aligned}\] 
It follows that $e^{1/2}f^{1/2}=f=f^{1/2}e^{1/2}$. Then 
\[(e^{1/2}-f^{1/2})^2=e-f. \]
Therefore $e\ge f$, and moreover $\widehat{E}_{\omega}(e-f)=0$. So $e=f$. 
\end{proof}

Propositions~\ref{P:equivalentconditions} and
\ref{P:equivalentcondition} show that the tracial Rokhlin property implies the
weak tracial Rokhlin property when $A$ is simple.

\begin{defn}[Hirshberg and Orovitz]
Let $G$ be a finite group and let $\alpha:G\to\Aut(A)$ be an action on a
simple, separable, unital $C\sp*$-algebra. The action $\alpha$ has the weak
tracial Rokhlin property if, for every nonzero positive element
$a\in A_{\omega}$, there are mutually orthogonal positive contractions
$\{e_g\}_{g\in G}\subset A_{\omega}\cap A'$ such that
\begin{enumerate}
\item $\alpha_{g, \omega}(e_h)=e_{gh}$ for every $g, h \in G$, 
\item $1-\sum_{g\in G} e_g \lesssim a $.
\end{enumerate}  
\end{defn}
\begin{ex}(\cite[Corollary 2.5]{GHS})
Let $A$ be a simple, finite, unital, infinite-dimensional $C\sp*$-algebra
with tracial rank zero and at most countably many extreme quasitraces, and let
$\alpha:G\to\Aut(A)$ be an action of a finite group. If
$\dim^c_{Rok}(\alpha)<\infty$, then $\alpha$ has the weak tracial Rokhlin
property. A concrete example is the cyclic-group action on a
higher-dimensional noncommutative torus constructed in
\cite[Proposition 2.8]{GHS}.
\end{ex}
\begin{ex}(\cite[Theorem 2.10]{GHS})
Let $A$ be a unital Kirchberg algebra and let $\alpha:G\to\Aut(A)$ be an
action of a finite group. Then $\alpha$ has the weak tracial Rokhlin property
if and only if it is pointwise outer; that is, $\alpha_g$ is not inner for
$g\ne1_G$.
\end{ex}
Given a finite-group action $\alpha:G\to\Aut(A)$, the natural conditional
expectation from $A$ onto the fixed-point algebra $A^{\alpha}$ is
\begin{equation}\label{Eq:average}
E(a)=\frac{1}{|G|} \sum_{g\in G} \alpha_g(a).
\end{equation}
If $\alpha$ has the weak tracial Rokhlin property, then it is pointwise outer;
the canonical expectation is therefore of index-finite type by
\cite{JP:saturated,W:index}.
\begin{thm}\label{T:action-inclusion}
Let $G$ be a finite group and let $\alpha:G\to\Aut(A)$ be an action on a
simple, separable, unital, infinite-dimensional $C\sp*$-algebra. Then the
following are equivalent:
\begin{enumerate}
\item $\alpha$ has the weak tracial Rokhlin property;
\item the conditional expectation $E:A\to A^{\alpha}$ in
\eqref{Eq:average} is of index-finite type and has the weak tracial Rokhlin
property.
\end{enumerate}
\end{thm}
\begin{proof}
Suppose that $\alpha$ has the weak tracial Rokhlin property. The preceding
observation shows that $E$ is of index-finite type. Given a nonzero positive
element $a\in A_{\omega}$, choose Rokhlin contractions
$\{e_g\}_{g\in G}$ and set $e=e_{1_G}$. Since $\Index E=|G|$,
\[\begin{aligned}
 (\Index E)E_{\omega}(e)e^{1/2}&=\sum_{g} \alpha_{g, \omega}(e_{1_G})e^{1/2}\\
&=(\sum_g  e_g) e^{1/2}\\
&=e^{3/2}.
\end{aligned}\] 
Moreover, \[1- (\Index E) E_{\omega}(e)=1- \sum_{g} e_g \lesssim a. \] 

Conversely, suppose that $E:A\to A^{\alpha}$ is of index-finite type and has
the weak tracial Rokhlin property. For a nonzero positive element
$a\in A_{\omega}$, choose a Rokhlin
contraction $e\in A_{\omega}\cap A'$ and set
$e_g=\alpha_{\omega,g}(e)$. Then
\[
 1-\sum_{g\in G}e_g=1-(\Index E)E_{\omega}(e)\lesssim a.
\]
Put $s=\sum_{g\in G}e_g$. Proposition~\ref{P:equivalentcondition} gives
$se^{1/2}=e^{3/2}$, and hence
\[
                    e^{1/2}(s-e)e^{1/2}=0.
\]
Since $s-e=\sum_{g\ne1_G}e_g$ is positive, it follows that
$e^{1/2}e_ge^{1/2}=0$, and hence $ee_g=0$, for every $g\ne1_G$. Applying
$\alpha_{\omega,h}$ shows that $e_he_g=0$ whenever $h\ne g$. Thus
$\{e_g\}_{g\in G}$ is a family of Rokhlin contractions for $\alpha$.
\end{proof}

Theorem~\ref{T:action-inclusion} covers fixed-point inclusions associated to
ordinary actions of finite groups. It is natural to ask whether there are
inclusions of index-finite type with the weak tracial Rokhlin property that
are not of this form. In Subsection~\ref{S:duality-applications}, we use the
duality theory to construct such inclusions; see
Proposition~\ref{P:nonabelian-dual-example}.

We next show that simplicity passes to the range of a weak tracial Rokhlin
expectation.
\begin{prop}\label{P:Rokhlinsimple}
Let $P\subset A$ be an inclusion of separable unital $C\sp*$-algebras of
index-finite type, with conditional expectation $E:A\to P$. If $A$ is simple
and $E$ has the weak tracial Rokhlin property, then $P$ is simple.
\end{prop}
\begin{proof}
First observe that a finite direct sum of unital simple $C\sp*$-algebras is
simple whenever its center is $\mathbb C$.

Let $B=C\sp*\langle A,e_P\rangle$ and let $\widehat E:B\to A$ be the dual
conditional expectation. Since $B$ is strongly Morita equivalent to $P$, it
suffices to prove that $B$ is simple. By Theorem~\ref{T:structure}, $B$ is a
finite direct sum of simple algebras, so it is enough to show that
$Z(B)=\mathbb C$.

Consider a Rokhlin positive contraction $e$.  Take $x\in B\cap A'$ of the form $\sum_i a_ie_Pb_i$. 
Then 
\[
\begin{split}
exe &= e(\sum_i a_ie_Pb_i)e\\
&=\sum a_iee_Peb_i\\
&=\sum_i a_i(\Index E)^{-1}e^2b_i\\
&=e \widehat{E}(x)e 
\end{split}
\] 

Since $\widehat E(x)\in A\cap A'$ for $x\in B\cap A'$, we obtain
\[ e(B\cap B')e \subset e(B\cap A')e \subset e(A\cap A')e= e \mathbb{C} e\]
We claim that the compression $x\mapsto exe$ is injective on $B\cap A'$.
Indeed, if $exe=0$, then
\[
 e\widehat E(x^*x)e=ex^*xe=x^*exe=0.
\]
The map $A\to A_{\omega}$ given by $a\mapsto eae$ is a nonzero order-zero
map and is therefore injective by Lemma~\ref{L:simpleorderzero}. Hence
$\widehat E(x^*x)=0$, and the faithfulness of $\widehat E$ gives $x=0$.
Since $eZ(B)e\subset\mathbb Ce^2$, this injectivity implies
$Z(B)=\mathbb C$. Therefore $B$, and hence $P$, is simple.
\end{proof}

Permanence from $A$ to $P$ is governed by a tracial approximate left inverse
for the inclusion $\iota:P\hookrightarrow A$. We now construct this inverse as
an order-zero map $A\to P_{\omega}$; compare
\cite[Proposition 4.22]{LO:2019}.

\begin{defn}\label{D:outer}
Let $1\in P\subset A$.  A conditional expectation $E:A\to P$ is
\emph{outer} if, for every $x\in A$ and every nonzero hereditary
$C\sp*$-subalgebra $C\subset P$,
\[
 \inf\{\|c(x-E(x))c\|:c\in C_+,\ \|c\|=1\}=0.
\]
\end{defn}
\begin{rem}\label{R:outer}
Suppose that $P\subset A$ is an irreducible inclusion of simple unital
$C\sp*$-algebras, that is, $P'\cap A=\mathbb C$, and that it has finite
depth.  Every conditional expectation $E:A\to P$ of index-finite type is outer by
\cite[Corollary 7.6]{Izumi:inclusion}.
\end{rem}

The outerness hypothesis appearing in the converse direction of the next
theorem is automatic for the two canonical expectations associated to an
outer finite-group action.  We record this fact separately for later use in
Proposition~\ref{P:nonabelian-dual-example}.

\newpage
\begin{lem}\label{L:outer-group-expectations}
Let $G$ be a finite group, let $A$ be a simple unital $C\sp*$-algebra, and
let $\alpha:G\to\Aut(A)$ be an outer action.  Then both canonical
conditional expectations
\[
 E_\alpha:A\longrightarrow A^\alpha,\qquad
 E_\alpha(a)=\frac1{|G|}\sum_{g\in G}\alpha_g(a),
\]
and
\[
 F_\alpha:A\rtimes_\alpha G\longrightarrow A,\qquad
 F_\alpha\left(\sum_{g\in G}a_gu_g\right)=a_{1_G},
\]
are outer.
\end{lem}

\begin{proof}
Pointwise outerness and simplicity of $A$ give
\[
 (A^\alpha)'\cap A=\mathbb C,
 \qquad
 A'\cap(A\rtimes_\alpha G)=\mathbb C;
\]
see \cite[Remark 3.13]{Izumi:inclusion}.  Thus both
\[
                       A^\alpha\subset A
 \quad\text{and}\quad
                       A\subset A\rtimes_\alpha G
\]
are irreducible.  The fixed-point algebra and the crossed product are simple
as well.  Moreover, pointwise outerness implies that the action is saturated,
so $A\rtimes_\alpha G$ is the basic construction for
$A^\alpha\subset A$; see \cite{SP:saturated}.  The shifted inclusion
\[
                         A\subset A\rtimes_\alpha G
\]
is of index-finite type and has depth two by
\cite[Lemma 3.1]{OT:tsr}.  Applying the Hilbert-bimodule identifications of
\cite[Remark 3.5]{JOPT:cancellation} to the shifted Jones tower, the boundedness
criterion in \cite[Proposition 3.10]{JOPT:cancellation} shows that the preceding
inclusion $A^\alpha\subset A$ has finite depth.  Thus both inclusions are of
index-finite type and have finite depth.  Applying
\cite[Corollary 7.6]{Izumi:inclusion} to the two inclusions shows that
$E_\alpha$ and $F_\alpha$ are outer.
\end{proof}

\begin{thm}\label{T:orderzerofromA}
Let $A$ be a finite simple unital $C\sp*$-algebra, let $P\subset A$ be a
unital inclusion, and let $E:A\to P$ be a conditional expectation of
index-finite type. If $E$ has the weak tracial Rokhlin property, then, for
every nonzero positive element $z\in P_{\omega}$, there are a nonzero positive
contraction $e\in A_{\omega}\cap A'$ and an injective c.p.c. order-zero map
$\beta:A\to P_{\omega}$ such that
\begin{enumerate}
\item $e=(\Index E) e^{1/2}e_Pe^{1/2}$,
\item $eae=\beta(a)e$ for all $a\in A$, 
\item $1-\beta(1) \lesssim z$ in $P_{\omega}$,
\item  $\beta(1)=g=(\Index E)E_{\omega}(e) \in P_{\omega}\cap P'$.
\end{enumerate}
Conversely, if $E$ is outer and such $e$ and $\beta$ exist for every nonzero
positive element of $P_{\omega}$, then $E$ has the weak tracial Rokhlin
property.
\end{thm}
\begin{proof}

First assume that $E$ has the weak tracial Rokhlin property. Regard
$z\in P_{\omega}$ as an element of $A_{\omega}$ and choose a Rokhlin
contraction $e$ such that, with $g=(\Index E)E_{\omega}(e)$,
$1-g\lesssim z$ in $A_{\omega}$. Define
$\beta(a)=(\Index E)E_{\omega}(ae)$. Since $e\in A_{\omega}\cap A'$,
$ae=e^{1/2}ae^{1/2}$, so $\beta$ is completely positive and contractive. For
$a\in A$,
\[\begin{split}
eae&= (\Index E) \widehat{E}_{\omega}(e_Pae^2) \\
&=(\Index E) \widehat{E}_{\omega}(e_Pae^{1/2}e e^{1/2}) \\
&=(\Index E) \widehat{E}_{\omega}(e_Pa e^{1/2}(\Index E)e^{1/2}e_Pe^{1/2}e^{1/2})  \quad \text{by Proposition~\ref{P:equivalentcondition}}\\
&=(\Index E) \widehat{E}_{\omega}(e_P a (\Index E) ee_Pe)\\
&=(\Index E)^2 \widehat{E}_{\omega}(E_{\omega}(ae)e_Pe)\\
&=(\Index E)E_{\omega}(ae)e. 
\end{split}\]
Let $a,b\in A_+$ be orthogonal. Then
\[\begin{split}
\beta(a)\beta(b)&=(\Index E)E_{\omega}(ae)(\Index E)E_{\omega}(be)\\
&=(\Index E)E_{\omega}((\Index E)E_{\omega}(ae)be)\\
&=(\Index E)E_{\omega}((\Index E)E_{\omega}(ae)eb)\\
&=(\Index E)E_{\omega}(eaeb)\\
&=(\Index E)E_{\omega}(eabe)=0.
\end{split}
\]
Similarly, $\beta(b)\beta(a)=0$, so $\beta$ has order zero. It is nonzero
and hence injective by Lemma~\ref{L:simpleorderzero}. The transfer of
$1-\beta(1)\lesssim z$ from $A_{\omega}$ to $P_{\omega}$ is precisely the
Cuntz-comparison argument in \cite[Theorem 4.16]{LeeOsaka4}. This proves
(1)--(4).

Conversely, assume that $E$ is outer and that the stated maps exist. Given a
nonzero positive element $z\in A_{\omega}$, the ultrapower version of
Lemma~\ref{L:contractioninsubalgebra} provides a nonzero positive element
$x\in P_{\omega}$ such that $x\lesssim z$ in $A_{\omega}$. Choose $e$ and
$\beta$ for $x$ and set $g=\beta(1)$. Then (1) gives the Rokhlin relation and
\[
                 1-g\lesssim x\lesssim z
                 \quad\text{in }A_{\omega}.
\]
Thus $e$ is a Rokhlin contraction for $z$.

 \end{proof}

\section{Weak tracial approximate representability for inclusions of index-finite type}
\label{S:approxrepresentability}
We introduce the notion dual to the weak tracial Rokhlin property. Recall that
a conditional expectation $E:A\to P$ is tracially approximately representable
if, for every nonzero positive element $z\in A_{\omega}$, there are projections
$e\in P_{\omega}\cap P'$ and $r\in A_{\omega}\cap A'$, together with a finite
set $\{u_i\}\subset A$, such that
\begin{enumerate}
\item $eae=E(a)e$ for all $a\in A$,
\item $\sum_i u_ieu_i^*=r$, and $re=e=er$, 
\item $1-r$ is Murray-von Neumann equivalent to a projection in $\overline{zA_{\omega}z}$ in $A_{\omega}$,
\item the map $P\ni x \mapsto xe$ is injective.
\end{enumerate}  
Replacing these projections by positive contractions leads to the following
definition.
\begin{defn}\label{D:approximaterepresentability}
Let $P\subset A$ be a unital inclusion and let $E:A\to P$ be a conditional
expectation of index-finite type. We say that $E$ is \emph{weakly tracially
approximately representable} if, for every nonzero positive element
$a\in A_{\omega}$, there are positive contractions
$e\in P_{\omega}\cap P'$ and $r\in A_{\omega}\cap A'$, together with a finite
set $\{u_i\}\subset A$, such that
\begin{enumerate}
\item $e^{1/2}xe^{1/2}=E(x)e$ for every $x\in A$, 
\item $\sum_i u_ieu_i^*=r$ and $re^{1/2}=e^{3/2}=e^{1/2}r$;
\item $1-r\lesssim a$ in $A_{\omega}$;
\item the map $x \mapsto xe$ is injective for $x\in P$.  
\end{enumerate}
\end{defn}

Weak tracial approximate representability is encoded by an order-zero map
from the basic construction. The following is the positive-contraction
analogue of \cite[Proposition 4.23]{LO:2019}.

\begin{prop}\label{P:approximaterepresentable}
Let $P\subset A$ be a unital inclusion, let $E:A\to P$ be of index-finite
type, and suppose that the basic construction
$B=C\sp*\langle A,e_P\rangle$ is simple. Then $E$ has the weak tracial
approximate representability if and only if, for every nonzero positive
element $z\in A_{\omega}$, there is an injective c.p.c. order-zero map
$\psi:B\to A_{\omega}$ such that
\begin{enumerate}
\item $\psi(a)=a\psi(1)=\psi(1)a$ for all $a\in A$ where $\psi(1) \in A_{\omega}\cap A'$, 
\item $\psi(e_P) \in P_{\omega} \cap P'$, 
\item $1-\psi(1) \lesssim z$.
\end{enumerate} 
\end{prop}
\begin{proof}
Let $B=C\sp*\langle A,e_P\rangle$. Suppose first that $E$ is weakly tracially
approximately representable. For a nonzero positive element
$z\in A_{\omega}$, choose positive contractions
$e\in P_{\omega}\cap P'$ and $r\in A_{\omega}\cap A'$, and a finite set
$\{u_i\}\subset A$, as in Definition~\ref{D:approximaterepresentability}.
Define $\psi(xe_Py)=xey$ for $x,y\in A$ and extend linearly. We first
verify that this formula is well defined. It is enough to show that
\[ \sum_i x_i e_P y_i=0 \Longrightarrow \sum_i x_i e y_i=0.\]
Since $\sum_i x_i e_P y_i=0$, 
\[\begin{split}
(\sum_i x_i e_P y_i )^* (\sum_j x_j e_P y_j )&=0\\
\sum_{i,j} y_i^*e_Px_i^*x_j e_P y_j&=0\\
\sum_{ij} y_i^*E(x_i^*x_j)e_Py_j&=0
\end{split}
\]  
Then \[\widehat{E}(\sum_{ij} y_i^*E(x_i^*x_j)e_Py_j)= (\Index E)^{-1}(\sum_{ij} y_i^*E(x_i^*x_j)y_j)=0.\]
It follows that $\sum_{i,j}y_i^*E(x_i^*x_j)y_j=0$ or 
\begin{equation}\label{E:welldefined}
\begin{pmatrix}
\begin{matrix}
y_1^* & \dots & y_n^* 
\end{matrix}\\
\bigzero
\end{pmatrix}
\begin{pmatrix}
E(x_1^*x_1) & \dots & E(x_1^*x_n) \\
\vdots & & \vdots\\
E(x_n^*x_1) &\dots & E(x_n^*x_n)
\end{pmatrix}
 \begin{pmatrix}
\begin{matrix}
y_1 \\
 \vdots \\
y_n 
\end{matrix} &
\bigzero
\end{pmatrix} =0.
\end{equation}
Note that $ \begin{pmatrix}
E(x_1^*x_1) & \dots & E(x_1^*x_n) \\
\vdots & & \vdots\\
E(x_n^*x_1) &\dots & E(x_n^*x_n)
\end{pmatrix}$ is positive and (\ref{E:welldefined}) implies that 
\[
\begin{pmatrix}
\begin{matrix}
y_1^* & \dots & y_n^* 
\end{matrix}\\
\bigzero
\end{pmatrix}\begin{pmatrix}
E(x_1^*x_1) & \dots & E(x_1^*x_n) \\
\vdots & & \vdots\\
E(x_n^*x_1) &\dots & E(x_n^*x_n)
\end{pmatrix}^{1/2}
 =0. \]
Thus 
\[\begin{split}
&(\sum_i x_i e y_i)^*(\sum_j x_j e y_j)\\
&=\sum_{i,j}y_i^*ex_i^*x_jey_j \\
&=\sum_{i,j} y_i^*E(x_i^*x_j)e^2y_j\\
&=\begin{pmatrix}
\begin{matrix}
y_1^* & \dots & y_n^* 
\end{matrix}\\
\bigzero
\end{pmatrix}
\begin{pmatrix}
E(x_1^*x_1) & \dots & E(x_1^*x_n) \\
\vdots & & \vdots\\
E(x_n^*x_1) &\dots & E(x_n^*x_n)
\end{pmatrix}
 \begin{pmatrix}
\begin{matrix}
e^2y_1 \\
 \vdots \\
e^2y_n 
\end{matrix} &
\bigzero
\end{pmatrix}=0
\end{split}\]
It follows that $\sum_i x_i ey_i=0$. \\
We next show that $\psi$ is completely positive. It suffices to consider a
Gram matrix $[z_i^*z_j]_{i,j}\in M_n(B)$. Write
$z_i=\sum_{k=1}^{n(i)}x_k^ie_Py_k^i$. Then
\[
 z_i^*z_j=
 \sum_{k=1}^{n(i)}\sum_{l=1}^{n(j)}
 (y_k^i)^*E((x_k^i)^*x_l^j)e_Py_l^j.
\]
If $w_i=\sum_{k=1}^{n(i)}x_k^ie^{1/2}y_k^i$, then
\[ \begin{split}
\psi(z_i^*z_j)&= \sum_{k=1}^{n(i)}\sum_{l=1}^{n(j)}(y_k^i)^*E((x_k^i)^*x_l^j)ey_l^j\\
&=\sum_{k=1}^{n(i)}\sum_{l=1}^{n(j)}(y_k^i)^*e^{1/2}(x_k^i)^*x_l^je^{1/2}y_l^j\\
&=w_i^*w_j
\end{split}
\]
This shows that $\psi$ is completely positive. 

Next we show that $\psi$ has order zero. Suppose that $u=\sum_{i=1}^nx_ie_Py_i \perp v=\sum_{j=1}^m r_je_Ps_j$. Then 
\[\begin{split}
uv=(\sum_{i=1}^nx_ie_Py_i)(\sum_{j=1}^m r_je_Ps_j)\\
=\sum_{i,j} x_iE(y_ir_j)e_Ps_j=0.
\end{split}\]
Thus 
\[\sum_{i,j}x_i E(y_ir_j)es_j=0.\]
Using $re=er=e^2$ and $r\in A_{\omega}\cap A'$, we have 
\[\begin{split}
\psi(u)\psi(v)&=(\sum_ix_i e y_i)(\sum_j r_jes_j)\\
&=\sum_{i,j}x_iE(y_ir_j)e^2s_j\\
&=\sum_{i,j}x_iE(y_ir_j)ers_j\\
&=(\sum_{i,j}x_iE(y_ir_j)es_j)r=0.
\end{split}\]
Similarly, $\psi(v)\psi(u)=0$, so $\psi$ has order zero. Since $B$ is
simple and $\psi$ is nonzero, Lemma~\ref{L:simpleorderzero} shows that
$\psi$ is injective.

We claim that the elements $u_i$ occurring in condition (2) form a
quasi-basis. First observe that, if $a\in A_+$ and $ae^2=0$, then $a=0$.
Since $\psi$ is an order-zero map, we can write $\psi(\cdot)=h\pi_{\psi}(\cdot)$. Then 
\[ae^2=0  \Longrightarrow \psi(ae_P)e=0.\]
By Lemma~\ref{L:lemma}, $\pi_{\psi}(E(a)e_P)e=0$. 
Thus, \begin{eqnarray*}
h\pi_{\psi}(E(a)e_P)e=0\\
\psi(E(a)e_P)e=0\\
E(a)e^2=0\\
eE(a)^{1/2}E(a)^{1/2}e=0\\
(E(a)^{1/2}e)^*(E(a)^{1/2}e)=0
\end{eqnarray*}
It follows that $E(a)e=0$. The injectivity of $P\ni x\mapsto xe$ and the
faithfulness of $E$ imply that $a=0$.

Using again $re=e^2=er$, note that for $a\in A$ 
\[\begin{split}
ae^2=are&=rae\\
&=\sum_i u_ieu_i^*ae\\
&=\sum_i u_iE(u_i^*a)e^2. 
\end{split}\]
It follows that $a=\sum_i u_iE(u_i^*a)$. Similarly,
$a=\sum_iE(au_i)u_i^*$. Hence $\{(u_i,u_i^*)\}$ is a quasi-basis, and
$\sum_i u_ie_Pu_i^*=1$. Therefore
\[
 \psi(1)=\sum_i u_i\psi(e_P)u_i^*=\sum_i u_ieu_i^*=r,
\]
so $1-\psi(1)=1-r\lesssim z$. This proves one direction.

Conversely, fix a nonzero positive element $z\in A_{\omega}$ and choose an
injective order-zero map
\[
              \psi:C\sp*\langle A,e_P\rangle\longrightarrow A_{\omega}
\]
satisfying conditions~(1)--(3) in the statement. Write
$\psi(\cdot)=\psi(1)\pi_{\psi}(\cdot)$, where $\pi_{\psi}$ is a
$*$-homomorphism, and set
\[
 e=\psi(e_P),\qquad r=\psi(1)\in A_{\omega}\cap A',
 \qquad e'=\pi_{\psi}(e_P).
\]
Then $e'$ is a projection.

Since $re'=e=e'r$, it follows that
$r^{1/2}e'=e^{1/2}=e'r^{1/2}$. Thus
\[
 re^{1/2}=r(r^{1/2}e')=r^{3/2}e'
 =r^{3/2}(e')^{3/2}=e^{3/2}.
\]
Similarly, $e^{1/2}r=e^{3/2}$. For a quasi-basis
$\{(u_i,u_i^*)\}$ for $E$, one has
\[
       \sum_i u_ie_Pu_i^*=1 \quad\text{in }C\sp*\langle A,e_P\rangle.
\]
Hence
\[\begin{split}
r=\psi(\sum_i u_i e_P u_i^*)&=\psi(1)\pi_{\psi}(\sum_i u_i e_P u_i^*)\\
&= \psi(1)\sum_i \pi_{\psi}(u_i) \pi_{\psi}(e_P)\pi_{\psi}(u_i^*) \\
&=\sum_i  \psi(u_i)\pi_{\psi}(e_P)\pi_{\psi}(u_i^*)\\
&=\sum_i u_i \pi_{\psi}(e_P) \psi(1) \pi_{\psi}(u_i^*)\\
&=\sum_i u_i \pi_{\psi}(e_P) \psi(1)u_i^*\\
&=\sum_i u_i \psi(e_P)u_i^*= \sum_i u_i e u_i^*.
\end{split}\]
The relation $e_Pxe_P=E(x)e_P$ gives
$e'\pi_{\psi}(x)e'=\pi_{\psi}(E(x))e'$ for $x\in A$. Consequently,
\[
\begin{split}
e^{1/2}xe^{1/2}
 &=e'r^{1/2}xr^{1/2}e'\\
 &=e'\psi(x)e'\\
 &=r e'\pi_{\psi}(x)e'\\
 &=r\pi_{\psi}(E(x))e'=\psi(E(x))e'=E(x)e.
\end{split}
\] 
In addition, $xe=0$ implies that $\psi(xe_P)=0$. Since $\psi$ is injective,
$xe_P=0$. Applying $\widehat E$ gives
$0=\widehat E(xe_P)=(\Index E)^{-1}x$, so $x=0$. Finally,
$1-r=1-\psi(1)\lesssim z$.

\end{proof}

\begin{rem}\label{R:quasibasis}
In the proof of Proposition~\ref{P:approximaterepresentable}, we showed that
the set $\{u_i\}$ gives rise to a quasi-basis for $E$.
\end{rem}

We recall tracial approximate representability in an ultrapower form. For
abelian groups this is \cite[Definition 4.8]{LO:2019}; for a general finite
group, the second condition below is replaced by its natural
conjugacy-covariant form.
\begin{defn}\label{D:tracialapproximaterepresentability}
Let $G$ be a finite group and let $A$ be a simple, separable, unital,
infinite-dimensional $C\sp*$-algebra. An action $\alpha:G\to\Aut(A)$ is
\emph{tracially approximately representable} if, for every nonzero positive
element $z\in A_{\omega}$, there are a projection
$e\in A_{\omega}\cap A'$ and a unitary representation
$w:G\to\mathcal U(eA_{\omega}e)$ such that
\begin{enumerate}
\item $\alpha_{\omega,g}(eae)=w_g(eae)w_g^*$ for all $a\in A$ and $g\in G$;
\item $\alpha_{\omega,g}(w_h)=w_{ghg^{-1}}$ for all $g,h\in G$;
\item $1-e$ is Murray--von Neumann equivalent to a projection in
$\overline{zA_{\omega}z}$.
\end{enumerate}
\end{defn}

\begin{thm}\cite[Corollary 4.25]{LO:2019}
Let $G$ be a finite abelian group, and let $\alpha$ be an outer action of
$G$ on a simple, separable, unital, infinite-dimensional $C\sp*$-algebra
$A$. Suppose that $A\rtimes_{\alpha}G$ is simple, and let $E$ be as in
Theorem~\ref{T:action-inclusion}. Then $\alpha$ is tracially approximately
representable if and only if $E$ is tracially approximately representable.
\end{thm}

For finite-group actions, we follow the local approximation picture introduced
by Asadi-Vasfi \cite[Definition 2.1]{AsadiVasfi}.

\begin{defn}\label{D:weaktracialapproxaction}
Let $G$ be a finite group, let $A$ be a simple, unital, infinite-dimensional
$C\sp*$-algebra, and let $\alpha:G\to\Aut(A)$ be an action. We say that
$\alpha$ is \emph{weakly tracially approximately representable} if, for every
finite set $\mathcal F\subset A$, every $\varepsilon>0$, and every
$x\in A_+$ with $\|x\|=1$, there are a positive contraction
$c\in A^\alpha$ and contractions $s_g\in\overline{cAc}$, for $g\in G$, such that
\begin{enumerate}
\item $\|s_{1_G}-c\|<\varepsilon$ and
$\|s_g^*-s_{g^{-1}}\|<\varepsilon$ for all $g\in G$;
\item $\|s_gs_h-cs_{gh}\|<\varepsilon$ for all $g,h\in G$;
\item $\|ca-ac\|<\varepsilon$ for every
$a\in\mathcal F\cup\{s_g:g\in G\}$;
\item $\|\alpha_g(cac)-s_gas_g^*\|<\varepsilon$ for every
$a\in\mathcal F$ and $g\in G$;
\item $\|\alpha_g(s_h)-s_{ghg^{-1}}\|<\varepsilon$ for all $g,h\in G$;
\item $(1-c-\varepsilon)_+\lesssim x$ in $A$;
\item $\|cxc\|>1-\varepsilon$.
\end{enumerate}
\end{defn}

The norm-largeness condition in Definition~\ref{D:weaktracialapproxaction} is
redundant in the stably finite case.

\begin{lem}\label{L:WTAR-norm-redundant}
Let $G$ be a finite group, let $A$ be a simple, unital, stably finite,
infinite-dimensional $C\sp*$-algebra with $T(A)\ne\varnothing$, and let
$\alpha:G\to\Aut(A)$ be an action.  Then $\alpha$ is weakly tracially
approximately representable if and only if, for every finite set
$\mathcal F\subset A$, every $\varepsilon>0$, and every $x\in A_+$ with
$\|x\|=1$, there are $c$ and $(s_g)_{g\in G}$ satisfying conditions
$(1)$--$(6)$ of Definition~\ref{D:weaktracialapproxaction}.
\end{lem}

\begin{proof}
Only the reverse implication requires proof.  We may assume that
$0<\varepsilon<1$.  By \cite[Lemma 2.9]{Phillips:large}, there is a positive
contraction
\[
             0\ne y\in\overline{xAx}
\]
such that, whenever $r\in A_+$ is a contraction with $r\lesssim y$, one has
\[
             \|(1-r)x(1-r)\|>1-\frac{\varepsilon}{3}.
\]
After rescaling, we may suppose that $\|y\|=1$.  Choose
$0<\delta<\varepsilon/3$ and apply the assumed six conditions to
$\mathcal F$, $\delta$, and $y$.  Let $c\in A^\alpha$ and
$(s_g)_{g\in G}$ be the resulting contractions, and put
\[
                    r=(1-c-\delta)_+.
\]
Then $r\lesssim y\lesssim x$.  Moreover, scalar functional calculus gives
\[
                 \|c-(1-r)\|\leq\delta.
\]
Consequently,
\[
\begin{aligned}
 \|cxc\|
 &\geq \|(1-r)x(1-r)\|-2\delta\\
 &>1-\frac{\varepsilon}{3}-2\delta
 >1-\varepsilon.
\end{aligned}
\]
Thus condition~(7) holds.  Conditions~(1)--(5) with tolerance $\varepsilon$
follow from their versions with tolerance $\delta$, and
\[
 (1-c-\varepsilon)_+
 \lesssim (1-c-\delta)_+
 =r\lesssim y\lesssim x,
\]
which gives condition~(6) for the originally prescribed element $x$.
\end{proof}

When $G$ is abelian, condition (5) reduces to
$\|\alpha_g(s_h)-s_h\|<\varepsilon$.  Notice also that approximate
representability implies Definition~\ref{D:weaktracialapproxaction} by taking
$c=1$.

We use a sequence-algebra consequence of this local definition only to connect
actions with inclusions of index-finite type.  Recall that a map $u:G\to B$ is an order
zero representation if each $u_g$ is a normal contraction, $u_{1_G}$ is
positive, and
\[
 u_gu_h=u_{1_G}u_{gh},
 \qquad
 u_g^*=u_{g^{-1}}
 \quad (g,h\in G).
\]
When $A$ is separable, a standard diagonal-sequence argument applied to
Definition~\ref{D:weaktracialapproxaction} gives, for every nonzero
$z\in A_\omega$, a positive contraction $e\in A_\omega\cap A'$ and an order
zero representation $u:G\to A_\omega$ such that
\[
 u_{1_G}=e,\qquad 1-e\lesssim z,
\]
\[
 \alpha_{\omega,g}(ea)=u_geau_g^*,\qquad
 \alpha_{\omega,g}(u_h)=u_{ghg^{-1}}.
\]
Indeed, one applies the local definition to increasing finite sets,
$\varepsilon_n\searrow0$, and normalized representatives of $z$; condition
(7) ensures that the resulting element $e$ is nonzero, while condition (6)
gives the Cuntz-small defect.  This is the sequence formulation used in
\cite[Definition 3.13]{Lee}.

Let $U:G\to\mathcal{U}(A\rtimes_{\alpha}G)$ denote the canonical implementing
unitary representation.  Thus $\alpha_g(a)=U_gaU_g^*$, and $\Ad U$ denotes
the action of $G$ on $A\rtimes_{\alpha}G$ given by conjugation by the
unitaries $U_g$.  The following characterization will be used below.

\begin{prop}\cite[Theorem 3.14]{Lee}\label{P:weaktracialapproxaction}
Let $G$ be a finite group, let $A$ be a unital separable $C\sp*$-algebra, and
let $\alpha:G\to\Aut(A)$ be an action.  Then the sequence conditions displayed
above hold if and only if, for every nonzero positive element
$z\in A_{\omega}$, there exists a nonzero equivariant completely positive
contractive order-zero map
\[
 \psi:(A\rtimes_{\alpha}G,\Ad U)\longrightarrow
       (A_{\omega},\alpha_{\omega})
\]
such that
\begin{enumerate}
\item $\psi(a)=a\psi(1)$ for every $a\in A$, with
$\psi(1)\in A_{\omega}\cap A'$;
\item $1-\psi(1)\lesssim z$ in $A_{\omega}$.
\end{enumerate}
\end{prop}

Combining Proposition~\ref{P:weaktracialapproxaction} with the order-zero
characterization of weak tracial approximate representability for conditional
expectations gives the following bridge from actions to inclusions.
  
\begin{thm}\cite[Corollary 3.16]{Lee}
\label{T:approximaterepresentableaction}
Let $G$ be a finite group, let $A$ be a simple, separable, unital,
infinite-dimensional $C\sp*$-algebra, and let
$\alpha:G\to\Aut(A)$ be an outer action. Let
\[
 E_{\alpha}:A\longrightarrow A^{\alpha},\qquad
 E_{\alpha}(a)=\frac{1}{|G|}\sum_{g\in G}\alpha_g(a),
\]
be the canonical conditional expectation.  If $\alpha$ is weakly tracially
approximately representable, then the conditional expectation $E_{\alpha}$ is
weakly tracially approximately representable.
\end{thm}

\begin{proof}
Set $P=A^{\alpha}$ and $B=A\rtimes_{\alpha}G$.  Since $\alpha$ is outer and
$A$ is simple, $B$ is simple.  Under the canonical identification
\[
 C\sp*\langle A,e_P\rangle\cong B,
\]
the Jones projection corresponds to
\[
 p=\frac{1}{|G|}\sum_{g\in G}U_g.
\]
Fix a nonzero positive element $z\in A_{\omega}$.  The diagonal-sequence
consequence of Definition~\ref{D:weaktracialapproxaction} and Proposition
\ref{P:weaktracialapproxaction} give a nonzero equivariant completely positive
contractive order-zero map
\[
 \psi:(B,\Ad U)\longrightarrow(A_{\omega},\alpha_{\omega})
\]
such that, with $r=\psi(1)$,
\[
 \psi(a)=ar=ra\quad(a\in A),\qquad
 r\in A_{\omega}\cap A',\qquad 1-r\lesssim z.
\]
The map $\psi$ is injective by Lemma~\ref{L:simpleorderzero}.  Write
$\psi(x)=r\pi_{\psi}(x)$, where $\pi_{\psi}$ is the supporting
$*$-homomorphism supplied by Theorem~\ref{T:orderzero}, and put
$f=\psi(p)$.  Since $(\Ad U)_g(p)=p$ for all $g\in G$, equivariance gives
$\alpha_{\omega,g}(f)=f$.  For a finite group,
$(A_{\omega})^{\alpha_{\omega}}=(A^{\alpha})_{\omega}$, so $f\in P_{\omega}$.
Moreover, if $a\in P$, then $pa=ap$, and
\[
\begin{aligned}
 fa
 &=\pi_{\psi}(p)ra
  =\pi_{\psi}(p)r\pi_{\psi}(a)
  =r\pi_{\psi}(pa) \\
 &=r\pi_{\psi}(ap)
  =r\pi_{\psi}(a)\pi_{\psi}(p)
  =af.
\end{aligned}
\]
Thus $f\in P_{\omega}\cap P'$.  Consequently, $\psi$ satisfies all three
conditions of Proposition~\ref{P:approximaterepresentable}, and therefore
$E_{\alpha}$ is weakly tracially approximately representable.
\end{proof}
 
\section{Dualities between Rokhlin property and approximate representability}\label{S:tracial}

\subsection{Duality theory}\label{S:duality-theory}

We first record an ultrapower lemma. It will be used in both duality
arguments, and every subsequent result states its own hypotheses explicitly.

\begin{lem}\label{L:ultrapowerorthogonalelements}
Let $D$ be a simple unital $C\sp*$-algebra which is not of type~I. For every
$a\in(D_{\omega})_+\setminus\{0\}$ and every positive integer $m$, there are
mutually orthogonal norm-one positive elements
\[
              r_1,\dots,r_m\in\overline{aD_{\omega}a}
\]
which are unitarily equivalent in $D_{\omega}$ and hence Cuntz equivalent.
\end{lem}
\begin{proof}
After scaling, represent $a$ by positive contractions $(a_k)_k$. Choose
$\delta>0$ such that $(a-2\delta)_+\ne0$. After changing the representatives
outside a set belonging to $\omega$, we may assume that
$c_k=(a_k-2\delta)_+$ is nonzero for every $k$.

The case $m=1$ is immediate, so assume that $m\ge2$. For each $k$, set
$H_k=\overline{c_kDc_k}$. This is a full hereditary subalgebra of $D$ and
hence is not of type~I. Apply \cite[Lemma 2.2]{Phillips:large}, with $m-1$ in
place of $n$, to the unital algebra
\[
                         \widetilde H_k=H_k+\mathbb C1_D.
\]
We obtain a unitary $u_k\in\widetilde H_k$ and a nonzero positive element
$b_k\in\widetilde H_k$ such that
\[
 b_k,\ u_kb_ku_k^*,\ \dots,\ u_k^{m-1}b_k(u_k^*)^{m-1}
\]
are mutually orthogonal. If $H_k\ne D$, apply the quotient map
$\widetilde H_k\to\widetilde H_k/H_k\cong\mathbb C$ to the orthogonality
relation $b_k(u_kb_ku_k^*)=0$. Since the quotient is commutative, the image
of $b_k$ has square zero and hence is zero. Thus $b_k\in H_k$; this is automatic
when $H_k=D$. After normalizing, we may assume that $\|b_k\|=1$.

Let
\[
             u=[(u_k)_k]\in D_{\omega},\qquad r_1=[(b_k)_k],
\]
and, for $1\le j\le m$, define
\[
                         r_j=u^{j-1}r_1(u^*)^{j-1}.
\]
Equivalently,
\[
             r_j=\big[(u_k^{j-1}b_k(u_k^*)^{j-1})_k\big].
\]
The elements $r_1,\dots,r_m$ are mutually orthogonal norm-one positive
elements.

Choose a continuous function $h$ on $[0,1]$ which vanishes near zero and is
one on $[2\delta,1]$. Since $c_k=(a_k-2\delta)_+$, functional calculus gives
\[
                         h(a_k)c_k=c_k.
\]
Consequently, $h(a_k)x=xh(a_k)=x$ for every
$x\in H_k=\overline{c_kDc_k}$. Each coordinate representing $r_j$ belongs
to $H_k$, because $H_k$ is an ideal in $\widetilde H_k$. It follows that
\[
                         r_j=h(a)r_jh(a).
\]
Finally, $h(a)\in\overline{aD_{\omega}a}$ because $h$ vanishes near zero.
Since $\overline{aD_{\omega}a}$ is hereditary, the displayed equality implies
that $r_j\in\overline{aD_{\omega}a}$ for every $j$.
\end{proof}
\begin{lem}\label{L:inclusiontechnical}
Let $P\subset A$ be an inclusion of unital $C\sp*$-algebras of index-finite
type, with conditional expectation $E:A\to P$. Let
$e\in A_{\omega}\cap A'$ be a positive contraction satisfying
\[
                   (\Index E)e^{1/2}e_Pe^{1/2}=e.
\]
If $p,q\in(P_{\omega})_+$ satisfy $pe=ep$ and
$q\lesssim e^2p$ in $A_{\omega}$, then $q\lesssim p$ in $P_{\omega}$.
\end{lem}
\begin{proof}
By assumption, there is a sequence $(v_k)_k$ in $A_{\omega}$ such that
\[
              \lim_{k\to\infty}\|v_ke^2pv_k^*-q\|=0.
\]
Set
\[
              w_k=\sqrt{\Index E}\,E_{\omega}(v_ke)\in P_{\omega}.
\]
Multiplying the Rokhlin relation on the left and right by $e^{1/2}$ gives
\[
                      (\Index E)ee_Pe=e^2.
\]
Using the Jones relation in the ultrapower of the basic construction, together
with $pe_P=e_Pp$ and $pe=ep$, we obtain
\[
\begin{aligned}
w_kpw_k^*e_P
 &= (\Index E)E_{\omega}(v_ke)pE_{\omega}(ev_k^*)e_P\\
 &= (\Index E)e_Pv_ke\,e_Pp\,ev_k^*e_P\\
 &= e_Pv_ke^2pv_k^*e_P\\
 &= E_{\omega}(v_ke^2pv_k^*)e_P.
\end{aligned}
\]
The map $P_{\omega}\ni z\mapsto ze_P$ is injective; indeed,
$\widehat E_{\omega}(ze_P)=(\Index E)^{-1}z$. Consequently,
\[
                  w_kpw_k^*=E_{\omega}(v_ke^2pv_k^*)\longrightarrow q.
\]
Thus $q\lesssim p$ in $P_{\omega}$.
\end{proof}

The next lemma is the positive-element analogue of
\cite[Lemma 4.21]{LO:2019}.

\begin{lem}\label{L:contractioninsubalgebra}
Let $P\subset A$ be an inclusion of unital $C\sp*$-algebras of index-finite
type, and let $E:A\to P$ be an outer conditional expectation. For every
$a\in A_+\setminus\{0\}$, there is $b\in P_+\setminus\{0\}$ such that
$b\lesssim a$.
\end{lem}
\begin{proof}

Fix $a\in A_+\setminus\{0\}$ and choose
$0<\varepsilon<\|E(a)\|/2$. By outerness, there is a positive contraction
$c\in P$ such that
\[
 \|c(a-E(a))c\|<\varepsilon,\qquad
 \|cE(a)c\|\ge\|E(a)\|-\varepsilon.
\]
See the proof of Theorem~2.1 in \cite{Osaka:(SP)-property} for details.
Set $b=(cE(a)c-\varepsilon)_+$. The norm estimate shows that $b\ne0$, and
\cite[Lemma 2.5]{KR:positive} gives
\[
                 b\lesssim cac\lesssim a.
\]
\end{proof}
\begin{lem}\label{L:normcontrol}
Let $A$ be a simple, finite, unital, infinite-dimensional $C\sp*$-algebra,
and let $x\in A_+$ satisfy $\|x\|=1$. For every $\varepsilon>0$, there is
$y\in\overline{xAx}_+\setminus\{0\}$ such that, whenever
$g\in A_{\omega}$ is a positive contraction with $g\lesssim y$ in
$A_{\omega}$, one has
\[
                       \|(1-g)x(1-g)\|>1-\varepsilon.
\]
\end{lem}
\begin{proof}
The ultrapower $A_{\omega}$ is finite: an isometry in $A_{\omega}$ lifts to
an approximate isometry and hence, by \cite[Lemma 6.2.4]{Ro:Classification},
to isometries in $A$. Use the elements denoted by $a,b$ and
$0\le d_j\le c_j\le b_j\le1$ in the proof of
\cite[Lemma 2.9]{Phillips:large}, and set $y=d_1$. If
$\|(1-g)x(1-g)\|\le1-\varepsilon$, the same argument gives
$1\lesssim (1-c_1-c_2)+d_1$ in $A_{\omega}$, contradicting finiteness.
\end{proof}
\raggedbottom
\begin{thm}\label{T:tracialduality1}
Let $P\subset A$ be an inclusion of unital, separable, infinite-dimensional
$C\sp*$-algebras, and let $E:A\to P$ be an outer conditional expectation of
index-finite type. Let $B=C\sp*\langle A,e_P\rangle$ and let
$\widehat E:B\to A$ be the dual conditional expectation. Suppose that $A$ is
simple and that $\widehat E$ is outer. Then $E$ has the weak tracial Rokhlin
property if and only if $\widehat E$ is weakly tracially approximately
representable.
\end{thm}
\begin{proof}
Assume that $E$ has the weak tracial Rokhlin property and fix
$z\in (B_{\omega})_+\setminus\{0\}$. Proposition~\ref{P:Rokhlinsimple}
implies that $P$ is simple, and hence $B$ is simple by
\cite[Corollary 2.2.14]{W:index}. Apply the ultrapower version of
Lemma~\ref{L:contractioninsubalgebra} first to $A\subset B$ and
$\widehat E$, and then to $P\subset A$ and $E$. We obtain a nonzero positive
element $p\in P_{\omega}$ such that $p\lesssim z$ in $B_{\omega}$.

The algebra $P$ is infinite-dimensional: otherwise the quasi-basis formula
would imply that $A$ is finite-dimensional. Let
$\{(v_i,v_i^*)\}_{i=1}^n$ be a quasi-basis for $E$. By
Lemma~\ref{L:ultrapowerorthogonalelements}, there are mutually orthogonal and
Cuntz equivalent positive contractions
$r_1,\dots,r_n\in\overline{pP_{\omega}p}$.
Choose a Rokhlin contraction $e\in A_{\omega}\cap A'$ such that
$1-(\Index E)E_{\omega}(e)\lesssim r_1$ in $P_{\omega}$. Since the $r_i$
are Cuntz equivalent, the same defect is Cuntz subequivalent to every $r_i$
in $P_{\omega}$.

Set $u_i=\sqrt{\Index E}\,v_ie_P$ and
$g=(\Index E)E_{\omega}(e)$. Then
$\displaystyle \sum_{i=1}^n u_i eu_i^*=\sum_i v_i g e_Pv_i^*.$
Thus, for any $x\in A$,
\[
\begin{aligned}
\sum_i u_i\widehat{E}(u_i^*x)e
 &=\sum_i (\Index E)v_ie_P\widehat{E}(e_Pv_i^*x)e\\
 &=\left(\sum_i v_i e_P v_i^*\right)xe=xe,\\
\sum_i e\widehat{E}(xu_i)u_i^*&=ex.
\end{aligned}
\]

Since $ge_P=e_Pg$ and $g\in P_{\omega}\cap P'$, define
$r=\sum_i u_ieu_i^*$. We show that $r\in B_{\omega}\cap B'$. Since
$B=C\sp*\langle A,e_P\rangle$ is generated by $A$ and $e_P$, it is enough to
check that $r$ commutes with both.
For $a\in A$,
\[
\begin{aligned}
ra&=\sum_i v_i ge_Pv_i^*a\\
&=  \sum_i v_i ge_P \left(\sum_k E(v_i^*a v_k)v_k^*\right)\\
&= \sum_i v_i\left(\sum_k E(v_i^*a v_k)ge_Pv_k^*\right)\\
&=\sum_k \sum_i v_i E(v_i^*av_k)ge_Pv_k^*\\
&=\sum_kav_kge_Pv_k^*=ar,\\
re_P&=\sum_i v_i ge_Pv_i^*e_P=  \sum_i v_i g E(v_i^*)e_P\\
&= \sum_i v_iE(v_i^*)e_Pg=ge_P,\\
e_Pr&=\sum_i e_Pv_i e_Pgv_i^*=  \sum_i E(v_i)e_Pgv_i^*\\
&= \sum_i ge_PE(v_i)v_i^*=ge_P.
\end{aligned}
\]

Note that $(\Index E)e^{1/2}e_Pe^{1/2}=e$.  Hence, for $x,y\in A$,
\[
\begin{aligned}
e^{1/2}(xe_Py)e^{1/2}
 &=xe^{1/2}e_Pe^{1/2}y
  =(\Index E)^{-1}xye \\
 &=\widehat{E}(xe_Py)e.
\end{aligned}
\]
By linearity and continuity, $e^{1/2}be^{1/2}=\widehat{E}(b)e$ for every
$b\in B$.
Let $\{e_{ij}\}_{i,j=1}^n$ be matrix units for $M_n$. In
$B_{\omega}\otimes M_n$, set
\[
V=\sum_i v_i\otimes e_{1i}
\quad\text{and}\quad
q=\sum_k(1-g)e_P\otimes e_{kk}.
\]
Then
all comparisons in the following display take place in
$B_{\omega}\otimes M_n$:
\[\begin{split}
(1-r)\otimes e_{11}&=VqV^*\\
&\sim q^{1/2}V^*Vq^{1/2}\\
&\lesssim q\\
&\lesssim \sum_k r_ke_P\otimes e_{kk}\\
&\sim (\sum_k r_k)e_P\otimes e_{11}\\
&\lesssim p^{1/2}e_Pp^{1/2}\otimes e_{11}\\
&\lesssim p\otimes e_{11}.
\end{split} \]
Hence we conclude that $1-r \lesssim p \lesssim z$ in $B_{\omega}$.

Conversely, suppose that $\widehat E$ has the weak tracial approximate
representability property. Fix $z\in(A_{\omega})_+\setminus\{0\}$. By the ultrapower
version of Lemma~\ref{L:contractioninsubalgebra}, choose a nonzero positive
contraction $p\in P_{\omega}$ with $p\lesssim z$ in $A_{\omega}$. Apply weak
tracial approximate representability to the nonzero element $pe_Pp\in
B_{\omega}$. We obtain positive contractions $e\in A_{\omega}\cap A'$ and
$r\in B_{\omega}\cap B'$, and a finite set $\{u_i\}\subset B$, such that
\begin{equation}\label{condition1}
e^{1/2}be^{1/2}=\widehat{E}(b)e \quad (b\in B),
\end{equation}
\begin{equation}\label{condition2}
\sum_i u_ieu_i^*=r,\qquad re^{1/2}=e^{3/2}=e^{1/2}r,
\end{equation}
\begin{equation}\label{condition3}
1-r  \lesssim pe_Pp \quad\text{in }B_{\omega}.
\end{equation}

From (\ref{condition1}) we have $(\Index E)e^{1/2}e_Pe^{1/2}=e$.  Set
\[
w_i=(\Index E)\widehat{E}(e_Pu_i)\in A.
\]
Then 
\[\begin{split}
e\left(\sum_i w_iw^*_i\right) e&= \sum_i (\Index E)^2\widehat{E}(e_Pu_i)e \widehat{E}(u_i^*e_P)e \quad (e\in A_{\omega}\cap A')\\
&=\sum_i (\Index E)^2 e^{1/2}e_P u_ieu_i^*e_Pe^{1/2}\\ 
&=(\Index E)^2 e^{1/2}e_P\left(\sum_i u_i eu_i^*\right)e_Pe^{1/2}\\
&=(\Index E)^2e^{1/2}e_Pre_Pe^{1/2}\\
&=(\Index E)^2e^{1/2}re_Pe^{1/2}\\
&=(\Index E)^2r e^{1/2}e_Pe^{1/2}\\
&=(\Index E)re=(\Index E) e^2. 
\end{split}
\]
Since $A$ is simple, the order-zero map $A\ni x\mapsto exe$ is injective.
It follows that $\sum_i w_iw_i^*=\Index E$.  Consequently,
\[
E_{\omega}\left(\sum_i w_iew_i^*\right)
 =E_{\omega}\left(\left[\sum_i w_iw_i^*\right]e\right)
 =(\Index E)E_{\omega}(e).
\]
We let $g=(\Index E)E_{\omega}(e)$.  Write
\[
u_i=\sum_j a_{ij}e_Pb_{ij}\in B,
\qquad
w_i=\sum_j E(a_{ij})b_{ij}.
\]
Then
\[ \begin{split}
ge_P=E_{\omega}(\sum_i w_iew_i^*)e_P&=e_P(\sum_i w_iew_i^*)e_P\\
&=\sum_i e_P\left(\sum_j E(a_{ij})b_{ij}\right)e \left(\sum_k b^*_{ik}E(a^*_{ik})\right)e_P\\
&=\sum_i \left(\sum_j e_Pa_{ij}e_Pb_{ij}\right)e\left(\sum_k b^*_{ik}e_Pa^*_{ik}e_P\right)\\
&=\sum_i e_Pu_ieu_i^*e_P\\
&=e_Pre_P=re_P.
\end{split}\]
 From the condition (\ref{condition3}), we have, in $B_{\omega}$,
\begin{eqnarray*}
(1-g)e_P=(1-r)e_P \lesssim pe_Pp.
\end{eqnarray*}
Therefore, for every $\varepsilon>0$, there is $q\in B_{\omega}$ such that
\begin{equation}\label{E:Cuntz}
 \|qpe_Ppq^* -(1-g)e_P\| < \varepsilon.
\end{equation}
Approximating $q$ by the algebraic span of $A_{\omega}e_PA_{\omega}$, we may
write $q=\sum_i a_ie_Pb_i$ with $a_i,b_i\in A_{\omega}$. Then
\eqref{E:Cuntz} becomes
\begin{equation}\label{E:Cuntzfinal}
\| \sum_i a_iE_{\omega}(b_ip)e_P \sum_k E_{\omega}(pb_k^*)a_k^* - (1-g)e_P\| <\varepsilon
\end{equation}
Taking $\widehat{E}_{\omega}$ in (\ref{E:Cuntzfinal}), we have 
\[\| \sum_i a_iE_{\omega}(b_i)p^2 \sum_k E_{\omega}(b_k^*)a_k^* - (1-g)\| <(\Index E)\varepsilon\]
This shows that $1-g \lesssim p^2$ in $A_{\omega}$. 
Consequently, in $A_{\omega}$ 
\[ 1-g \lesssim p^2 \lesssim p \lesssim z.\]
\end{proof}
\flushbottom

The converse duality requires the additional assumption that $A$ is finite.
  
\begin{thm}\label{T:tracialduality2}
Let $P\subset A$, $E:A\to P$, and
$B=C\sp*\langle A,e_P\rangle$ satisfy the assumptions of
Theorem~\ref{T:tracialduality1}. Assume in addition that $A$ is finite. Then
$E$ is weakly tracially approximately representable if and only if
$\widehat E:B\to A$ has the weak tracial Rokhlin property.
\end{thm}
\begin{proof}
Suppose that $E$ is weakly tracially approximately representable. Fix
$z\in(B_{\omega})_+\setminus\{0\}$. By the ultrapower version of
Lemma~\ref{L:contractioninsubalgebra}, choose a nonzero positive contraction
$p\in A_{\omega}$ with $p\lesssim z$ in $B_{\omega}$. Apply weak tracial approximate
representability to $p$ to obtain a positive contraction
$e\in P_{\omega}\cap P'$ satisfying
\begin{equation}\label{E:approxcondition(1)}
e^{1/2}xe^{1/2}=E(x)e \quad \text{for all $x \in A$},
\end{equation}
a positive contraction $r\in A_{\omega}\cap A'$ and a finite set
$\{u_i\}\subset A$ such that
\begin{equation}\label{E:approxcondition(2)}
\sum_i u_i eu_i^*=r,   re^{1/2}=e^{3/2}=e^{1/2}r,
\end{equation} 
\begin{equation}\label{E:approxcondition(3)}
1-r \lesssim p \quad\text{in }A_{\omega}.
\end{equation}
Note that $\{(u_i,u_i^*)\}$ is a quasi-basis for $E$ by Remark~\ref{R:quasibasis}. 
Define $f$ in $B_{\omega}$ by 
\[ f=\sum_i u_ie e_P u_i^*.\] 
Note that
\[ \begin{split}
fr&=\sum_i  u_i ee_Pu_i^* \sum_j u_jeu_j^*\\
&=\sum_{i,j} u_ie_Peu_i^*u_jeu_j^*\\
&=\sum_j \left(\sum_i u_i E(u_i^* u_j)e^2\right)e_Pu_j^* \quad \text{by (\ref{E:approxcondition(1)})}\\
&=\sum_j u_je^2  e_Pu_j^* =f^2 \quad \text{$\{(u_i,u^*_i)\}$ a quasi-basis for $E$}.
\end{split}\]
Thus $f^{1/2}r^{1/2}=f=r^{1/2}f^{1/2}$. On the other hand, for any  $a\in A$
\[ \begin{split}
fa&=\sum_i  u_i ee_Pu_i^*a\\
&=\sum_i u_ie_Peu_i^*a  \quad (e_Pee_P=E_{\omega}(e)e_P=ee_P)\\
&=\sum_i u_ie_P \left(\sum_j eE(u_i^*a u_j)u_j^*\right) \quad  \text{$\{(u_i,u^*_i)\}$ a quasi-basis for $E$}\\
&=\sum_j \left(\sum_i u_iE(u_i^*au_j)e \right) e_Pu_j^* \quad \text{by $e \in P_{\omega}\cap P'$}\\
&=\sum_j au_j ee_Pu_j^*=af \quad  \text{$\{(u_i,u^*_i)\}$ a quasi-basis for $E$}.
\end{split}\]
Moreover, 
\[ \begin{split}
fe_P&=\left(\sum_i  u_i ee_Pu_i^*\right)e_P\\
&=\sum_i u_ieE(u_i^*)e_P\\
&=\sum_i \left(u_i E(u_i^*)e\right)e_P \\
&=ee_P \quad  \text{$\{(u_i,u^*_i)\}$ a quasi-basis for $E$}\\
\end{split}
\]
and 
\[\begin{split}
e_Pf&=e_P \left(\sum_i u_i e_Pe u_i^* \right)\\
&=\sum_i e_Pu_ie_P e u_i^*\\
&=\sum_i E(u_i)e_Peu_i^*\\
&=\sum_i e_Pe E(u_i)u_i^*\\
&=e_Pe \quad  \text{$\{(u_i,u^*_i)\}$ a quasi-basis for $E$}.
\end{split}\]
Thus $f\in B_{\omega}\cap B'$. Moreover,
\[\widehat{E}_{\omega}(f)=(\Index E)^{-1}\left(\sum_i u_i eu_i^*\right).\]
Therefore $(\Index E)\widehat{E}_{\omega}(f)=r$ is a positive contraction such that
 \[\begin{split}
f^{1/2}(\Index E)\widehat{E}_{\omega}(f)&=f^{1/2}r^{1/2}r^{1/2}\\
&=fr^{1/2}=f^{1/2}f^{1/2}r^{1/2}\\
&=f^{1/2}f=f^{3/2},
\end{split}\] and  
 \[1-r\lesssim p\lesssim z\quad\text{in }B_{\omega}.\]
Thus $f$ is a Rokhlin contraction for $\widehat E$ by
Proposition~\ref{P:equivalentcondition}.
 
Conversely, suppose that $\widehat E$ has the weak tracial Rokhlin property
and fix $z\in(A_{\omega})_+\setminus\{0\}$. Since $A$ is simple and not of
type~I, Lemma~\ref{L:ultrapowerorthogonalelements} provides orthogonal norm-one
positive elements $r_1,r_2\in\overline{zA_{\omega}z}$ with $r_1\sim r_2$.

A standard diagonal-sequence argument, as in
\cite[Theorem 3.13]{LeeOsaka}, gives a positive contraction
$f'\in B_{\omega}\cap B'$ such that
\[
 f'r_1=r_1f',\qquad
 (\Index E)\widehat E_{\omega}(f')(f')^{1/2}=(f')^{3/2},
 \qquad (f')^2r_1\ne0.
\]
For completeness, represent $r_1$ by norm-one positive elements
$r_1^{(k)}\in A$. Apply Lemma~\ref{L:normcontrol}, with a fixed number
$0<\varepsilon<1$, inside
$\overline{(r_1^{(k)})^2A(r_1^{(k)})^2}$, and choose a Rokhlin contraction
whose defect is Cuntz subequivalent to the resulting element. The norm
estimate in Lemma~\ref{L:normcontrol}, together with the Rokhlin relation,
gives a positive lower bound for the norm of the product of this contraction
with $r_1^{(k)}$; the bound is independent of $k$. Diagonalizing over an
increasing sequence of finite subsets of $B$ gives $f'$ with the three
displayed properties. This is the diagonalization used in
\cite[Theorem 3.13]{LeeOsaka}.

Now we invoke a Rokhlin positive contraction $f\in B_{\omega}\cap B'$ such that
 $1-(\Index E)\widehat{E}_{\omega}(f) \lesssim (f')^2r_1$ in $B_{\omega}$. 

Define $e=(\Index E)\widehat{E}_{\omega}(fe_P)\in A_{\omega}$. Since
$fe_P=e_Pf$, we may write $fe_P=pe_P$ for some $p\in P_{\omega}$. Then
\[ e=(\Index E)\widehat{E}_{\omega}(fe_P)=(\Index E)\widehat{E}_{\omega}(pe_P)=p. \] 
Since $e$ commutes with any element in $P$ by the definition, $e\in P_{\omega}\cap P'$.

Note that \[
\begin{split}
&(f^{1/2}e_P)(f^{1/2}e_P)=fe_P \\
&=ee_P=(e^{1/2}e_P)(e^{1/2}e_P).
\end{split}
\]
It follows that $(fe_P)^{1/2}=f^{1/2}e_P=e^{1/2}e_P$ and $e^{1/2}=(\Index E)\widehat{E}_{\omega}(f^{1/2}e_P)$.

For $a\in A$
\[ 
\begin{split}
e^{1/2} a e^{1/2}&=(\Index E)^2 \widehat{E}_{\omega}(f^{1/2}e_P)a \widehat{E}_{\omega}(f^{1/2}e_P)\\
&=(\Index E)^2 \widehat{E}_{\omega}(f^{1/2}e_P)\widehat{E}_{\omega}(af^{1/2}e_P)\\
&=(\Index E)\widehat{E}_{\omega}(f^{1/2}e_P (\Index E)\widehat{E}_{\omega}(af^{1/2}e_P))\\
&=(\Index E)\widehat{E}_{\omega}(e_Pe^{1/2}e_P (\Index E)\widehat{E}_{\omega}(af^{1/2}e_P))\\
&=(\Index E)\widehat{E}_{\omega}(e_Pe^{1/2} (\Index E)\widehat{E}_{\omega}(af^{1/2}e_P)e_P )\\
&=(\Index E)\widehat{E}_{\omega}(e_Pe^{1/2} a(\Index E)\widehat{E}_{\omega}(f^{1/2}e_P)e_P)\\
&=(\Index E)\widehat{E}_{\omega}(e^{1/2}e_Pae_Pe^{1/2})\\
&=(\Index E)\widehat{E}_{\omega}(e^{1/2}E(a)e_Pe^{1/2})\\
&=E(a)e.
\end{split}
\]
  Then take a quasi-basis $\{(u_i,u_i^*)\}$ for $E$. Note that $\sum_i u_i e_P u^*_i =1$.  It follows that \[
\begin{split}
\sum_i u_i e u_i^*&=(\Index E)\sum_i u_i \widehat{E}_{\omega}(fe_P)u_i^*\\
&=(\Index E)\widehat{E}_{\omega}(\sum_i u_ife_Pu_i^*)\\
&=(\Index E)\widehat{E}_{\omega}(f). 
\end{split}
\]   
Thus $r=\sum_i u_i e u_i^*$ is in $A_{\omega}\cap A'$. Moreover,  
\[\begin{split}
re^{1/2}&=\sum_i u_ieu_i^*e^{1/2}\\
&=\sum_i u_i e^{1/2}(e^{1/2}u_i^*e^{1/2})\\
&=\sum_i u_i e^{1/2}E(u_i^*)e\\
&=\sum_i u_iE(u_i^*)e^{3/2}\\
&=e^{3/2}.
\end{split}
\]
Similarly, $e^{1/2}r=e^{3/2}$. 
Apply Lemma~\ref{L:inclusiontechnical} to the inclusion $A\subset B$, the
elements $1-r$ and $r_1$, and the contraction $f'$. Since
$1-r\lesssim(f')^2r_1$ in $B_{\omega}$ and $f'$ commutes with $r_1$, we get
\[
                  1-r\lesssim r_1\sim r_2\lesssim z
                  \quad\text{in }A_{\omega}.
\]
It remains to verify the injectivity condition in
Definition~\ref{D:approximaterepresentability}. Since $B$ is simple and
$0\ne f\in B_{\omega}\cap B'$, the order-zero map $b\mapsto bf$ is injective.
If $x\in P$ and $xe=0$, then
\[
                 (xe_P)f=xf e_P=xe e_P=0.
\]
Hence $xe_P=0$, and applying $\widehat E$ gives $x=0$. Therefore
$E:A\to P$ is weakly tracially approximately representable.
\end{proof}
\subsection{Applications: inclusions beyond ordinary group actions}
\label{S:duality-applications}

The non-group-type inclusions in \cite{LOT:Rokhlin-index} arise from quantum
symmetries on $\mathcal O_2$ and $\mathcal O_\infty$ and are detected by their
non-integer index.  The following proposition gives a complementary
integer-index mechanism in the stably finite setting: duality converts a
non-abelian weakly tracially approximately representable action into an
inclusion of index-finite type whose associated dual Hopf algebra is not a
group algebra, and hence which is not a fixed-point inclusion arising from an
ordinary finite-group action.

We recall the Hopf-algebraic picture behind this mechanism.  Let \(\mathcal H\)
be a finite-dimensional Hopf \(C\sp*\)-algebra acting on a unital
\(C\sp*\)-algebra \(M\), let \(h_{\mathcal H}\) be its normalized Haar
idempotent, and put
\[
 P=M^{\mathcal H},\qquad
 E_{\mathcal H}(x)=h_{\mathcal H}\mathbin{\cdot}x \quad (x\in M).
\]
Szyma\'nski and Peligrad proved that, when the action is saturated,
\(E_{\mathcal H}:M\to P\) is of index-finite type with
\(\Index E_{\mathcal H}=\dim(\mathcal H)1_M\), and that the Haar idempotent
\(h_{\mathcal H}\in M\rtimes\mathcal H\) satisfies the Jones relations for
\(E_{\mathcal H}\); see \cite{SP:saturated}.  Watatani's universal property
of the basic construction \cite[Chapter II]{W:index} therefore gives a
$*$-homomorphism
\[
 \Theta:C\sp*\langle M,e_P\rangle\longrightarrow M\rtimes\mathcal H,
 \qquad \Theta(xe_Py)=xh_{\mathcal H}y.
\]
Saturation says that the closed linear span of
\(Mh_{\mathcal H}M\) is the whole crossed product, so \(\Theta\) is
surjective.  Moreover, the canonical conditional expectation
\(F_{\mathcal H}:M\rtimes\mathcal H\to M\) satisfies, on spanning elements,
\[
 F_{\mathcal H}(\Theta(xe_Py))
   =F_{\mathcal H}(xh_{\mathcal H}y)
   =\frac{1}{\dim(\mathcal H)}xy
   =\widehat E_{\mathcal H}(xe_Py).
\]
Thus \(F_{\mathcal H}\circ\Theta=\widehat E_{\mathcal H}\).  Since the dual
conditional expectation \(\widehat E_{\mathcal H}\) is faithful, \(\Theta\)
is injective and hence is an isomorphism.  The equality with the dual
conditional expectation is therefore a direct calculation, not an additional
appeal to the crossed-product theorem.
The dual Hopf \(C\sp*\)-algebra
\(\mathcal H^0\) acts on \(M\rtimes\mathcal H\), and its fixed-point algebra
is \(M\); see \cite{SP:saturated}.

In particular, an action \(\alpha:G\curvearrowright A\) of a finite group is
equivalently an action of the group Hopf \(C\sp*\)-algebra
\[
 \mathcal H=C\sp*(G)=\operatorname{span}\{U_g:g\in G\},
 \qquad \Delta(U_g)=U_g\otimes U_g,
\]
and
\[
 A\rtimes\mathcal H\cong A\rtimes_\alpha G.
\]
Here and below, \(C\sp*(G)\) denotes the group algebra, whereas \(C(G)\)
denotes its dual Hopf \(C\sp*\)-algebra of functions on \(G\).  Thus, in the
saturated case, the tower
\[
                 A^\alpha\subset A\subset A\rtimes_\alpha G
\]
is the basic-construction tower, and the coefficient expectation
\[
 F\left(\sum_{g\in G}a_gu_g\right)=a_{1_G}
\]
is both the Haar expectation for the dual \(C(G)\)-action and the dual
conditional expectation of \(E_\alpha\).

\begin{rem}\label{R:examples-after-duality}
The hypothesis is stated at the weakest of the three representability levels
considered here.  Since approximate representability implies tracial
approximate representability, which in turn implies weak tracial approximate
representability, the proposition applies uniformly at all three levels.  This
allows the same duality mechanism to convert strict separations among the
three action properties into examples of inclusions.
\end{rem}

\begin{prop}\label{P:nonabelian-dual-example}
Let $G$ be a finite non-abelian group, and let
$\alpha:G\curvearrowright A$ be an outer weakly tracially approximately
representable action on a simple, separable, unital, finite,
infinite-dimensional $C\sp*$-algebra $A$. Let
\[
 F:A\rtimes_{\alpha}G\longrightarrow A
\]
be the canonical conditional expectation, given by
\[
 F\left(\sum_{g\in G}a_gu_g\right)=a_{1_G}.
\]
Then the inclusion
\[
                  A\subset A\rtimes_{\alpha}G
\]
has the weak tracial Rokhlin property with respect to $F$.

Moreover, this inclusion is not of ordinary group type, even up to an
isomorphism preserving the specified conditional expectations. More precisely,
there do not exist a finite group $H$, an action
$\beta:H\curvearrowright D$, and a $*$-isomorphism
\[
 \Phi:A\rtimes_{\alpha}G\longrightarrow D
\]
such that
\[
 \Phi(A)=D^{\beta}
 \quad\text{and}\quad
 E_{\beta}\circ\Phi=\left.\Phi\right|_A\circ F,
\]
where
\[
 E_{\beta}(d)=\frac{1}{|H|}\sum_{h\in H}\beta_h(d),
 \qquad d\in D.
\]
\end{prop}

\begin{proof}
Let
\[
 E_{\alpha}:A\longrightarrow A^{\alpha},
 \qquad
 E_{\alpha}(a)=\frac{1}{|G|}\sum_{g\in G}\alpha_g(a),
\]
be the canonical conditional expectation. The basic construction for the
inclusion $A^{\alpha}\subset A$ is canonically isomorphic to the crossed
product:
\[
 C\sp*\langle A,e_{A^{\alpha}}\rangle\cong A\rtimes_{\alpha}G.
\]
Indeed, the group Hopf $C\sp*$-algebra $C\sp*(G)$ acts on $A$ through
$\alpha$, and the action is saturated because it is outer.  The
Szyma\'nski--Peligrad identification recalled above therefore applies, and
under it the dual conditional expectation $\widehat{E}_{\alpha}$ corresponds
to $F$.  By
Lemma~\ref{L:outer-group-expectations}, both $E_\alpha$ and $F$ are outer,
so the outerness hypotheses in Theorem~\ref{T:tracialduality2} are
satisfied. Since $\alpha$ is weakly tracially approximately representable,
Theorem~\ref{T:approximaterepresentableaction} shows directly that
$E_{\alpha}$ is weakly tracially approximately representable. Therefore
Theorem~\ref{T:tracialduality2} implies that $F$ has the weak tracial
Rokhlin property.

It remains to prove the second assertion. The finite-dimensional Hopf
$C\sp*$-algebra canonically associated with the dual inclusion
\[
 A\subset A\rtimes_{\alpha}G
\]
is the dual Hopf $C\sp*$-algebra $(C(G),\Delta_G)$, where
\[
 \Delta_G(f)(g_1,g_2)=f(g_1g_2),
 \qquad f\in C(G),\quad g_1,g_2\in G.
\]
Suppose, toward a contradiction, that $H$, $\beta$, and $\Phi$ exist as in
the statement. Replacing $H$ by $H/\ker(\beta)$ changes neither $D^\beta$
nor $E_\beta$, so we may assume that $\beta$ is faithful. The inclusion
$A\subset A\rtimes_\alpha G$ is irreducible by
\cite[Remark 3.13]{Izumi:inclusion}; hence $D^\beta\subset D$ is irreducible
as well. If $\beta_h=\Ad(v)$ for a unitary $v\in D$, then $v$ commutes with
$D^\beta$ and is therefore scalar. Faithfulness of $\beta$ then gives
$h=1_H$. Thus $\beta$ is outer and hence saturated. Moreover, $D$ and
$D^\beta$ are simple because they are isomorphic to
$A\rtimes_\alpha G$ and $A$, respectively.

By \cite[Lemma 3.1]{OT:tsr}, the inclusion
$A\subset A\rtimes_\alpha G$ has depth two. Since $\Phi$ preserves the
specified conditional expectations, it extends to the first basic
constructions by
\[
 \widetilde\Phi(xe_Fy)=\Phi(x)e_{E_\beta}\Phi(y),
 \qquad x,y\in A\rtimes_\alpha G,
\]
and then recursively to an isomorphism of the two basic-construction towers.
Thus $D^\beta\subset D$ also has depth two.
The reconstruction theorem \cite[Theorem 3.1]{Szlachanyi} says that an
inclusion of index-finite type and depth two with finite-dimensional centers
determines a normalized regular weak Hopf $C\sp*$-algebra, uniquely. In the
present irreducible setting, the source and target base algebras are scalar,
so the reconstructed objects are ordinary Hopf $C\sp*$-algebras. Here they
are $C(G)$ for $A\subset A\rtimes_\alpha G$ and $C\sp*(H)$ for
$D^\beta\subset D$. The tower isomorphism would therefore induce a Hopf
$C\sp*$-algebra isomorphism
\[
 C\sp*(H)\cong C(G).
\]
Every group Hopf $C\sp*$-algebra $C\sp*(H)$ is cocommutative.  On the other
hand, the coproduct on $C(G)$ is cocommutative if and only if $G$ is
abelian.  Since $G$ is non-abelian, $C(G)$ is not cocommutative and hence
cannot be isomorphic to $C\sp*(H)$ as a Hopf $C\sp*$-algebra.
Thus the inclusion $A\subset A\rtimes_{\alpha}G$, equipped with $F$, is not
isomorphic to a fixed-point inclusion arising from an action of an ordinary
finite group.
\end{proof}

We now use Proposition~\ref{P:nonabelian-dual-example} as a conversion
principle.  In each of the next three constructions, we first produce a
pointwise outer action of a finite non-abelian group at a different level of
representability: an approximately representable action on an infinite-type
UHF algebra, a tracially approximately representable but not approximately
representable action on a simple monotracial AF algebra, and a weakly
tracially approximately representable action on $\mc Z$ which is neither
tracially nor approximately representable.  The corollary following each
construction then applies Proposition~\ref{P:nonabelian-dual-example} to turn
the action into a non-group-type inclusion of index-finite type with the weak
tracial Rokhlin property.  We begin with the strongest of the three
representability properties.

\begin{ex}[A non-abelian AR action on an infinite-type UHF algebra]
\label{E:UHF-AR}
Let $G$ be a finite non-abelian group, put $m=|G|$, and let
\[
             \lambda:G\longrightarrow\mathcal U(\ell^2(G))
\]
be the left regular representation.  Using
$B(\ell^2(G))\cong M_m(\mathbb C)$, set
\[
 A=\bigotimes_{n=1}^{\infty}M_m(\mathbb C)\cong M_{m^\infty},
 \qquad
 \alpha_g=\bigotimes_{n=1}^{\infty}\Ad(\lambda_g).
\]
Then $A$ is an infinite-type UHF algebra and $\alpha:G\curvearrowright A$
is approximately representable and pointwise outer.

Indeed, for $N\in\mathbb N$ define
\[
 v_g^{(N)}=\lambda_g^{\otimes N}\otimes1
       \in\mathcal U(A).
\]
If $a$ belongs to the first $N_0$ tensor factors, then, for $N\geq N_0$,
\[
       \Ad(v_g^{(N)})(a)=\alpha_g(a).
\]
Moreover, for all $g,h\in G$ and every $N$,
\[
 v_g^{(N)}v_h^{(N)}=v_{gh}^{(N)},
 \qquad
 \alpha_g(v_h^{(N)})=v_{ghg^{-1}}^{(N)}.
\]
Consequently the unitaries
\[
             v_g=[(v_g^{(N)})_N]\in\mathcal U(A_\omega)
\]
implement $\alpha_g$, form a unitary representation of $G$, and satisfy the
conjugacy-covariance relation
$\alpha_{\omega,g}(v_h)=v_{ghg^{-1}}$.  These are precisely the
sequence-algebra relations for approximate representability.

It remains to verify pointwise outerness.  We give explicitly the norm-Cauchy
obstruction underlying the product-type criterion
\cite[Theorem 6.3]{EvansKawahigashi}.  Fix $g\ne1_G$, let
$d=\operatorname{ord}(g)$, and put
\[
             V_N=\lambda_g^{\otimes N}\otimes1.
\]
Write $A_N=M_m^{\otimes N}$ and let $E_N:A\to A_N$ be the canonical
trace-preserving conditional expectation.  Suppose that
$\alpha_g=\Ad(U)$ for a unitary $U\in A$.  Since
$\Ad(U)|_{A_N}=\Ad(V_N)|_{A_N}$, we have
\[
                  V_N^*U\in A_N'\cap A.
\]
Consequently, with $\mu_N=\tau(V_N^*U)\in\mathbb C$,
\[
                  E_N(U)=\mu_NV_N.
\]
The expectations $E_N$ converge pointwise in norm to the identity, and hence
$\mu_NV_N\to U$ in norm.  It follows that $|\mu_N|\to1$.  After replacing
$\mu_N$ by its phase $z_N=\mu_N/|\mu_N|$ for all sufficiently large $N$, we
obtain
\[
                         z_NV_N\longrightarrow U.
\]
In particular, $(z_NV_N)_N$ would be a norm-Cauchy sequence.

This is impossible.  Indeed, $V_{N+1}=V_N\otimes\lambda_g$, so
\begin{align*}
 \|z_{N+1}V_{N+1}-z_NV_N\|
   &=\|z_{N+1}\lambda_g-z_N1_{M_m}\|\\
   &=\|\lambda_g-\overline{z_{N+1}}z_N1_{M_m}\|\\
   &\geq \delta_g,
\end{align*}
where
\[
       \delta_g=\min_{z\in\mathbb T}\|\lambda_g-z1_{M_m}\|>0.
\]
The last inequality follows from the spectral decomposition: as a
permutation matrix, $\lambda_g$ is the direct sum of $m/d$ cycles of length
$d$, so its spectrum contains all $d$th roots of unity and is not a singleton.
Thus $\lambda_g$ is not a scalar unitary, and compactness of $\mathbb T$ gives
$\delta_g>0$.  This contradiction proves that $\alpha_g$ is outer.  Hence the
action is pointwise outer.
\end{ex}

\begin{cor}\label{C:UHF-nongroup-inclusion}
Let $G$ be a finite non-abelian group and let
$\alpha:G\curvearrowright M_{|G|^\infty}$ be the action in
Example~\ref{E:UHF-AR}.  Then the inclusion
\[
 M_{|G|^\infty}\subset M_{|G|^\infty}\rtimes_\alpha G,
\]
equipped with its canonical conditional expectation, has the weak tracial
Rokhlin property and is not isomorphic to a fixed-point inclusion arising
from an action of an ordinary finite group.
\end{cor}

\begin{proof}
Approximate representability implies weak tracial approximate
representability, and the action is pointwise outer.  The conclusion follows
from Proposition~\ref{P:nonabelian-dual-example}.
\end{proof}

\begin{ex}[A non-abelian TAR action which is not AR]
\label{E:nonabelian-TAR-not-AR}
Let $G$ be a finite non-abelian group and put $m=|G|$.  We construct a simple
unital AF algebra $A$ and an action $\alpha:G\curvearrowright A$ which is
tracially approximately representable but not approximately representable.

Choose positive integers $k_1,r_1$ and a sequence of positive integers
$\ell_n\to\infty$.  Inductively set
\[
 k_{n+1}=\ell_nk_n+mr_n,
 \qquad
 r_{n+1}=k_n+r_n,
\]
and let
\[
 A_n=B_n\oplus D_n,
 \qquad
 B_n=M_{k_n}(\mathbb C),
 \qquad
 D_n=\bigoplus_{h\in G}M_{r_n}(\mathbb C).
\]
Let $U^{(n)}:G\to\mathcal U(\mathbb C^{k_n})$ and
$W^{(n)}:G\to\mathcal U(\mathbb C^{r_n})$ be unitary representations.  On
$B_n$ and $D_n$, respectively, define
\[
 \alpha_g^{(n)}(x)=U_g^{(n)}xU_g^{(n)*},
 \qquad
 \bigl(\alpha_g^{(n)}(y)\bigr)_k
 =W_g^{(n)}y_{g^{-1}k}W_g^{(n)*}.
\]
Writing $\lambda$ for the left regular representation of $G$, choose the
representations at the next stage by
\[
 U_g^{(n+1)}
 =\bigl(U_g^{(n)}\otimes1_{\ell_n}\bigr)
  \oplus\bigl(\lambda_g\otimes W_g^{(n)}\bigr),
 \qquad
 W_g^{(n+1)}=U_g^{(n)}\oplus W_g^{(n)}.
\]
For $x\in B_n$ and $y=(y_h)_{h\in G}\in D_n$, set
\begin{equation}\label{E:AF-connecting-map}
 \phi_n(x,y)=
 \left(
  \operatorname{diag}\bigl(
   \underbrace{x,\ldots,x}_{\ell_n\text{ copies}},(y_h)_{h\in G}
  \bigr),
  \bigl(\operatorname{diag}(x,y_h)\bigr)_{h\in G}
 \right).
\end{equation}
The maps $\phi_n$ are unital, injective, and equivariant.  Every summand of
$A_n$ maps into $B_{n+1}$, and $B_{n+1}$ maps into every summand of
$A_{n+2}$.  Hence every two-step incidence matrix is strictly positive, so
\[
                 A=\varinjlim(A_n,\phi_n)
\]
is a simple unital AF algebra and the actions $\alpha^{(n)}$ induce an action
$\alpha:G\curvearrowright A$.

We next verify uniform concentration of traces on the main block.  Let
$p_n=1_{D_n}$, $e_n=1_{B_n}=1-p_n$, and $t_n=k_n/r_n$.  For every
$\tau\in T(A)$,
\begin{equation}\label{E:AF-trace-bound}
 \tau(p_n)\leq
 \max\left\{
   \frac{mr_n}{k_{n+1}},\frac{r_n}{r_{n+1}}
 \right\}.
\end{equation}
Moreover,
\[
 t_{n+1}=\frac{\ell_nt_n+m}{t_n+1}.
\]
Once $\ell_n\geq m$, one has $t_{n+1}\geq m$, and at subsequent stages
\[
 t_{n+1}\geq\frac{m}{m+1}\ell_n.
\]
Thus $t_n\to\infty$, and both terms on the right side of
\eqref{E:AF-trace-bound} tend to zero.  Consequently,
\[
                  \sup_{\tau\in T(A)}\tau(p_n)\longrightarrow0.
\]
Since $B_n$ is a single matrix summand, the possible discrepancy between two
traces on $A_n$ is supported on $D_n$.  The preceding uniform estimate shows
that $A$ has a unique tracial state, denoted by $\tau$.

We claim that $\alpha$ is tracially approximately representable in the
non-abelian sense of Definition~\ref{D:tracialapproximaterepresentability}.
Let $\mathcal F\subset A$ be finite, let $\varepsilon>0$, and let
$a\in A_+$ have norm one.  Approximate $\mathcal F\cup\{a\}$ inside some
$A_j$.  Since $A$ is a simple infinite-dimensional AF algebra, it has
Property~(SP); choose a nonzero projection $q\in\overline{aAa}$.  For
sufficiently large $N>j$,
\[
                         \tau(p_N)<\tau(q).
\]
Strict comparison of projections in $A$ gives $p_N\lesssim q$.  Set
\[
                 e=e_N,
 \qquad
                 u_g=U_g^{(N)}\in\mathcal U(eAe).
\]
The projection $e$ is fixed by $\alpha$, is central on $A_N$, and therefore
approximately commutes with $\mathcal F$.  On the main corner,
\[
 \alpha_g(exe)=u_gexe u_g^*,
 \qquad
 u_gu_h=u_{gh},
 \qquad
 \alpha_g(u_h)=u_{ghg^{-1}}.
\]
Furthermore, the compression of $A_j$ into $B_N$ is faithful, because the
main block eventually contains a copy of every summand of $A_j$.  Hence
$\|eae\|>1-\varepsilon$ after making the initial approximation sufficiently
close.

We spell out the passage to the ultrapower formulation.  Let
$0\ne z=[(z_n)_n]\in(A_{\omega})_+$, with the representatives chosen so that
their norms are bounded away from zero, and let $(\mathcal F_n)_n$ be an
increasing sequence of finite subsets with dense union in $A$.  Apply the
preceding construction coordinatewise to $\mathcal F_n$, $1/n$, and
$z_n/\|z_n\|$.  We obtain stages $N(n)$, projections
$e^{(n)}=e_{N(n)}$, and representations
$u^{(n)}:G\to\mathcal U(e^{(n)}Ae^{(n)})$, with $1-e^{(n)}$ equivalent to a
projection in $\overline{z_nAz_n}$.  Then
\[
 e=[(e^{(n)})_n],\qquad w_g=[(u_g^{(n)})_n]\qquad(g\in G)
\]
satisfy all three conditions of
Definition~\ref{D:tracialapproximaterepresentability}.  Thus $\alpha$ is
tracially approximately representable.

It remains to rule out approximate representability.  With minimal-rank
generators for the matrix summands,
\[
 K_0(A_n)=\mathbb Z\oplus\mathbb Z[G].
\]
Let $\varepsilon:\mathbb Z[G]\to\mathbb Z$ be the augmentation map,
$\Omega=\sum_{h\in G}h$, and $I_G=\ker(\varepsilon)$.  The connecting map is
\begin{equation}\label{E:AF-K0-map}
 (\phi_n)_*(x,y)=\bigl(\ell_nx+\varepsilon(y),x\Omega+y\bigr).
\end{equation}
In particular, $(\phi_n)_*(0,y)=(0,y)$ for $y\in I_G$, so $I_G$ embeds
canonically in $K_0(A)$.  The induced action is
\[
 (\alpha_g^{(n)})_*(x,y)=(x,\lambda_gy),
\]
and left translation by every $g\ne1_G$ acts nontrivially on $I_G$.  Thus
$(\alpha_g)_*\ne\id_{K_0(A)}$ for $g\ne1_G$.  Since every approximately
inner automorphism acts trivially on $K$-theory, $\alpha$ is pointwise outer
and is not approximately representable.
\end{ex}

\begin{cor}\label{C:AF-nongroup-inclusion}
Let $G$ be a finite non-abelian group and let
$\alpha:G\curvearrowright A$ be the action in
Example~\ref{E:nonabelian-TAR-not-AR}.  Then the inclusion
\[
                   A\subset A\rtimes_\alpha G,
\]
equipped with its canonical conditional expectation, has the weak tracial
Rokhlin property and is not isomorphic to a fixed-point inclusion arising
from an action of an ordinary finite group.
\end{cor}

\begin{proof}
Tracial approximate representability implies weak tracial approximate
representability.  The conclusion therefore follows from
Proposition~\ref{P:nonabelian-dual-example}.
\end{proof}

\begin{lem}\label{L:projectionless-ultrapower}
Let $D$ be a unital $C\sp*$-algebra whose only projections are $0$ and $1$.
Then the only projections in $D_{\omega}$ are $0$ and $1$.
\end{lem}

\begin{proof}
Let $p\in D_{\omega}$ be a projection.  Represent $p$ by self-adjoint
contractions $(x_n)_n$ with
$\lim_{n\to\omega}\|x_n^2-x_n\|=0$.  Functional calculus gives projections
$q_n\in D$ such that $\lim_{n\to\omega}\|x_n-q_n\|=0$.  Each $q_n$ is either
$0$ or $1$, and the ultrafilter contains exactly one of the sets on which
$q_n=0$ or $q_n=1$.  Hence $p=0$ or $p=1$.
\end{proof}

\begin{ex}[A WTAR action on $\mc Z$ which is neither TAR nor AR]
\label{E:Z-WTAR-not-TAR-not-AR}
Let $G$ be a finite non-abelian group and put $m=|G|$.  We construct a
pointwise outer action
\[
                         \alpha:G\curvearrowright\mc Z
\]
which is weakly tracially approximately representable but is neither
tracially nor approximately representable.

For relatively prime integers $p,q>1$, write
\[
 Z_{p,q}=
 \{f\in C([0,1],M_{pq}):
   f(0)\in M_p\otimes1_q,\ 
   f(1)\in1_p\otimes M_q\}.
\]
Choose distinct primes $p_1,q_1>m$, and write
\[
 p_1=a_1m+r_p,\qquad q_1=b_1m+r_q,
 \qquad 0<r_p,r_q<m.
\]
If $\lambda:G\to\mathcal U(\mathbb C^m)$ is the left regular
representation, set
\[
 U_g^{(1)}=\lambda_g^{\oplus a_1}\oplus1_{r_p},
 \qquad
 W_g^{(1)}=\lambda_g^{\oplus b_1}\oplus1_{r_q},
 \qquad
 X_g^{(1)}=U_g^{(1)}\otimes W_g^{(1)}.
\]
Thus $X_g^{(1)}$ is non-scalar whenever $g\ne1_G$.

Assume that $p_n$ and $q_n$ have been chosen and are relatively prime.
Dirichlet's theorem gives distinct primes
\[
 P_n=k_nq_n+1,\qquad Q_n=\ell_np_n+1.
\]
We choose them sufficiently large that
\begin{equation}\label{E:Z-trace-contraction}
                 \eta_n:=\frac1{P_n}+\frac1{Q_n}<2^{-n}.
\end{equation}
Set
\[
 p_{n+1}=p_nP_n,\qquad q_{n+1}=q_nQ_n,
 \qquad K_n=P_nQ_n.
\]
Then $\gcd(p_{n+1},q_{n+1})=1$ and
\[
 M_{p_{n+1}q_{n+1}}
 \cong M_{p_nq_n}\otimes M_{K_n}.
\]
Let $r_0$ and $r_1$ be the remainders of $K_n$ modulo $q_{n+1}$ and
$p_{n+1}$, respectively.  Since
\[
 K_n=k_nq_{n+1}+Q_n
     =\ell_np_{n+1}+P_n,
\]
we have
\begin{equation}\label{E:Z-remainders}
                         r_0=Q_n,\qquad r_1=P_n.
\end{equation}
For $1\le j\le K_n$, define
\begin{equation}\label{E:Z-eigenvalue-maps}
 \xi_{n,j}(t)=
 \begin{cases}
  t/2,     &1\le j\le r_0,\\
  1/2,     &r_0<j\le K_n-r_1,\\
  (t+1)/2, &K_n-r_1<j\le K_n.
 \end{cases}
\end{equation}

We record the endpoint multiplicities.  If
$f(0)=a_f\otimes1_{q_n}$, then
\[
 a_0:=\frac{q_nr_0}{q_{n+1}}=1,\qquad
 b_0:=\frac{K_n-r_0}{q_{n+1}}=k_n.
\]
Consequently there is a permutation unitary $u_n(0)$, independent of $f$,
such that
\begin{equation}\label{E:Z-left-endpoint}
\begin{split}
 &u_n(0)
 \left(\bigoplus_{j=1}^{K_n}f(\xi_{n,j}(0))\right)u_n(0)^*\\
 &\hspace{18mm}
 =\left(a_f^{\oplus a_0}\oplus
 f(1/2)^{\oplus b_0}\right)\otimes1_{q_{n+1}}
 \in M_{p_{n+1}}\otimes1_{q_{n+1}}.
\end{split}
\end{equation}
Indeed,
$p_na_0+p_nq_nb_0=p_n(1+k_nq_n)=p_{n+1}$.
Likewise, if $f(1)=1_{p_n}\otimes b_f$, then
\[
 a_1:=\frac{K_n-r_1}{p_{n+1}}=\ell_n,\qquad
 b_1:=\frac{p_nr_1}{p_{n+1}}=1,
\]
and a permutation unitary $u_n(1)$ satisfies
\begin{equation}\label{E:Z-right-endpoint}
\begin{split}
 &u_n(1)
 \left(\bigoplus_{j=1}^{K_n}f(\xi_{n,j}(1))\right)u_n(1)^*\\
 &\hspace{18mm}
 =1_{p_{n+1}}\otimes
 \left(f(1/2)^{\oplus a_1}\oplus b_f^{\oplus b_1}\right)
 \in1_{p_{n+1}}\otimes M_{q_{n+1}}.
\end{split}
\end{equation}
Here
$p_nq_na_1+q_nb_1=q_n(\ell_np_n+1)=q_{n+1}$.

Define the representations recursively by
\[
\begin{aligned}
 U_g^{(n+1)}&=U_g^{(n)}\otimes1_{P_n},&
 W_g^{(n+1)}&=W_g^{(n)}\otimes1_{Q_n},\\
 X_g^{(n+1)}&=U_g^{(n+1)}\otimes W_g^{(n+1)}.
\end{aligned}
\]
In the block-diagonal picture put
$Y_g=X_g^{(n)}\otimes1_{K_n}$.  The endpoint permutations can be chosen
equivariantly:
\begin{equation}\label{E:Z-endpoint-intertwiners}
                  X_g^{(n+1)}u_n(i)=u_n(i)Y_g
                  \quad(g\in G,\ i=0,1).
\end{equation}
For completeness, at the left endpoint the representation on the
regrouped $p_{n+1}$-dimensional block is
\[
 U_g^{(n)}\oplus
 (U_g^{(n)}\otimes1_{q_n})^{\oplus k_n}
 \cong U_g^{(n)}\otimes1_{1+k_nq_n}
 =U_g^{(n+1)},
\]
and its $q_{n+1}$-dimensional multiplicity space carries
$W_g^{(n)}\otimes1_{Q_n}=W_g^{(n+1)}$.  At the right endpoint the
$p_{n+1}$-dimensional multiplicity space carries $U_g^{(n+1)}$, while
the representation on the regrouped $q_{n+1}$-dimensional block is
\[
 (1_{p_n}\otimes W_g^{(n)})^{\oplus\ell_n}\oplus W_g^{(n)}
 \cong W_g^{(n)}\otimes1_{\ell_np_n+1}
 =W_g^{(n+1)}.
\]

It follows from \eqref{E:Z-endpoint-intertwiners} that
$u_n(0)^*u_n(1)$ belongs to the finite-dimensional commutant
$\{Y_g:g\in G\}'$.  Its unitary group is path connected, so choose
\[
 z_n:[0,1]\longrightarrow\mathcal U(\{Y_g:g\in G\}')
\]
with $z_n(0)=1$ and $z_n(1)=u_n(0)^*u_n(1)$, and put
$v_n(t)=u_n(0)z_n(t)$.  Then
\begin{equation}\label{E:Z-path-intertwines}
 X_g^{(n+1)}v_n(t)
 =v_n(t)(X_g^{(n)}\otimes1_{K_n}).
\end{equation}
Therefore
\begin{equation}\label{E:Z-connecting-map}
 \phi_n(f)(t)=v_n(t)
 \left(\bigoplus_{j=1}^{K_n}f(\xi_{n,j}(t))\right)v_n(t)^*
\end{equation}
defines a unital injective homomorphism
$Z_{p_n,q_n}\to Z_{p_{n+1},q_{n+1}}$.

We briefly verify that the limit is $\mc Z$.  The two nonconstant branches
in \eqref{E:Z-eigenvalue-maps} cover $[0,1]$, and their $d$-fold
compositions have image diameter at most $2^{-d}$ and meet every nonempty
open interval at all sufficiently late stages.  Thus the limit is simple.
For any two limit traces $\rho$ and $\sigma$, the common constant branch
$1/2$ disappears when their restrictions are subtracted, while the total
weight of the two remaining branches is $\eta_n$.  Hence, for
$a\in Z_{p_n,q_n}$ and $N>n$,
\[
 |\rho(a)-\sigma(a)|
 \le2\|a\|\prod_{j=n}^{N-1}\eta_j.
\]
By \eqref{E:Z-trace-contraction}, the right hand side tends to zero, so
the limit is monotracial.  The Jiang--Su inductive-limit criterion now
identifies
\[
                    \varinjlim(Z_{p_n,q_n},\phi_n)\cong\mc Z;
\]
see \cite[Propositions 2.5 and 2.8 and Theorem 6.2]{JiangSu}.

Let $\alpha_g^{(n)}=\Ad(X_g^{(n)})$.  These actions preserve the endpoint
fibers, and \eqref{E:Z-path-intertwines} gives
\[
 \alpha_g^{(n+1)}\circ\phi_n
 =\phi_n\circ\alpha_g^{(n)}.
\]
They therefore induce an action $\alpha:G\curvearrowright\mc Z$.

We next prove weak tracial approximate representability.  Let
$\mathcal F\subset\mc Z$ be finite, let $\varepsilon>0$, and let
$x\in\mc Z_+$ have norm one.  Approximate $\mathcal F\cup\{x\}$ at some
stage and pass to a stage $N$ such that
\[
              \eta_N<\min\{\tau(x),\varepsilon/3\},
\]
where $\tau$ is the unique trace of $\mc Z$.  If $\mu_N$ represents the
restriction of $\tau$ to $Z_{p_N,q_N}$, then
\[
 \tau(a)=\int_0^1\operatorname{tr}_{p_Nq_N}(a(t))\,d\mu_N(t).
\]
For $0<\delta<1/2$, \eqref{E:Z-eigenvalue-maps} and
\eqref{E:Z-remainders} give
\begin{equation}\label{E:Z-endpoint-trace}
 \mu_N([0,\delta]\cup[1-\delta,1])
 \le\frac{r_0+r_1}{K_N}
 =\frac1{P_N}+\frac1{Q_N}
 =\eta_N.
\end{equation}
Let $c_\delta$ be linear from $0$ to $1$ on $[0,\delta]$, equal to $1$
on $[\delta,1-\delta]$, and linear from $1$ to $0$ on
$[1-\delta,1]$.  Set
\[
 c=c_\delta1_{p_Nq_N},\qquad s_g=cX_g^{(N)}.
\]
The element $c$ is fixed by $\alpha$, and $s_g$ belongs to
$\overline{cZ_{p_N,q_N}c}$ because it vanishes at both endpoints.  On the
stage algebra the following relations hold exactly:
\[
 s_{1_G}=c,\qquad s_g^*=s_{g^{-1}},\qquad
 s_gs_h=cs_{gh},
\]
\[
 [c,s_g]=0,\qquad
 \alpha_g(cac)=s_gas_g^*,\qquad
 \alpha_g(s_h)=s_{ghg^{-1}}.
\]
They therefore hold within $\varepsilon$ on $\mathcal F$.  Moreover,
\eqref{E:Z-endpoint-trace} gives
\[
 d_\tau(1-c)\le\eta_N<\tau(x)\le d_\tau(x).
\]
Strict comparison in $\mc Z$ yields $1-c\lesssim x$, which is stronger
than condition~(6) of Definition~\ref{D:weaktracialapproxaction}.
Lemma~\ref{L:WTAR-norm-redundant} supplies the norm-largeness condition.
Thus $\alpha$ is weakly tracially approximately representable.

It remains to separate the stronger properties.  Since $\mc Z$ has no
nontrivial projections, Lemma~\ref{L:projectionless-ultrapower} shows that
neither does $\mc Z_\omega$.  If $\alpha$ were approximately representable, an implementing
representation $u:G\to\mathcal U(\mc Z_\omega)$ would satisfy
$u_g^k=1$ whenever $g$ has order $k$.  The spectral projections of $u_g$
would therefore lie in $\mc Z_\omega$, forcing $u_g$ to be scalar.
This would make every $\alpha_g$ trivial.  In contrast, for $g\ne1_G$,
the non-scalar matrix $X_g^{(1)}$ acts nontrivially on an element supported
in the interior of $Z_{p_1,q_1}$, and equivariance and injectivity preserve
this in the limit.  Hence $\alpha$ is not approximately representable.

The same argument shows that $\alpha$ is pointwise outer.  Indeed, if
$\alpha_g=\Ad(v)$ and $g$ has order $k$, then $v^k$ is scalar.  After
rescaling, $v^k=1$, and the absence of nontrivial spectral projections in
$\mc Z$ forces $v$ to be scalar, a contradiction.

Finally, choose a nonzero $z\in(\mc Z_\omega)_+$ with
$d_{\tau_\omega}(z)<1$.  In the ultrapower formulation of tracial
approximate representability, the supporting projection cannot be zero,
since otherwise $1$ would be Murray--von Neumann equivalent to a projection
in $\overline{z\mc Z_\omega z}$.  It must therefore be $1$, because
$\mc Z_\omega$ has no nontrivial projections.  Tracial approximate
representability would then reduce to approximate representability, already
ruled out.  Thus $\alpha$ is neither TAR nor AR.
\end{ex}

\begin{cor}\label{C:Z-nongroup-inclusion}
Let $G$ be a finite non-abelian group and let
$\alpha:G\curvearrowright\mc Z$ be the action in
Example~\ref{E:Z-WTAR-not-TAR-not-AR}.  Then
\[
                       \mc Z\subset\mc Z\rtimes_\alpha G,
\]
equipped with its canonical conditional expectation, has the weak tracial
Rokhlin property and is not isomorphic, as an inclusion with specified
conditional expectation, to a fixed-point inclusion arising from an action
of any ordinary finite group.
\end{cor}

\begin{proof}
The action is pointwise outer and weakly tracially approximately
representable, so the conclusion follows from
Proposition~\ref{P:nonabelian-dual-example}.
\end{proof}

\section{Applications via tracially sequentially-split order-zero maps}
\label{S:applications}

We now place all of the permanence results in the unified framework of
tracially sequentially-split maps by order zero. This framework was developed
in \cite{LeeOsaka4}; its consequences for tracial comparison, divisibility,
and tracial nuclear dimension were established in \cite{LeePermanenceII}.

\begin{defn}\cite[Definition 3.1]{LeeOsaka4}\label{D:sequentiallysplitbyorderzero}
Let $C$ and $D$ be unital $C\sp*$-algebras. A unital $*$-homomorphism
$\varphi:C\to D$ is said to be \emph{tracially sequentially-split by order
zero} if, for every nonzero positive element $z\in C_{\omega}$, there are a
c.p.c. order-zero map $\psi:D\to C_{\omega}$ and a nonzero positive
contraction $g\in C_{\omega}\cap C'$ such that
\[
       \psi(\varphi(c))=cg \quad(c\in C),
       \qquad 1-g\lesssim z \quad\text{in }C_{\omega}.
\]
The map $\psi$ is called a tracial approximate left inverse for $\varphi$.
\end{defn}

The order-zero map obtained from a Rokhlin contraction in
Theorem~\ref{T:orderzerofromA} is precisely such a tracial approximate left
inverse.

\begin{prop}\label{P:inclusion-sequentiallysplit}
Let $P\subset A$ be an inclusion of unital, separable, infinite-dimensional
$C\sp*$-algebras of index-finite type. Suppose that $A$ is simple and finite
and that the conditional expectation $E:A\to P$ has the weak tracial Rokhlin
property. Then the canonical inclusion $\iota:P\hookrightarrow A$ is
tracially sequentially-split by order zero.
\end{prop}

\begin{proof}
Fix a nonzero positive element $z\in P_{\omega}$. By
Theorem~\ref{T:orderzerofromA}, there is an injective c.p.c. order-zero map
\[
             \beta:A\longrightarrow P_{\omega},
             \qquad \beta(a)=(\Index E)E_{\omega}(ae),
\]
where $e\in A_{\omega}\cap A'$ is a Rokhlin contraction,
$g=\beta(1)\in P_{\omega}\cap P'$, and $1-g\lesssim z$ in $P_{\omega}$.
For $p\in P$, the $P$-bimodule property of $E_{\omega}$ gives
\[
                  \beta(\iota(p))=p\beta(1)=pg.
\]
Thus $\beta$ is a tracial approximate left inverse for $\iota$ in the sense
of Definition~\ref{D:sequentiallysplitbyorderzero}.
\end{proof}

For the terminology in the next theorem, we use tracial $m$-comparison and
tracial $m$-almost divisibility in the sense of Winter, and tracial nuclear
dimension in the sense recalled in \cite[Definition 3.7]{LeePermanenceII}.

\begin{thm}\label{T:framework-permanence}
Let $P\subset A$ and $E:A\to P$ satisfy the assumptions of
Proposition~\ref{P:inclusion-sequentiallysplit}. Then the following hold.
\begin{enumerate}
\item If $A$ is tracially $\mc{Z}$-absorbing, then $P$ is tracially
$\mc{Z}$-absorbing.
\item If $A$ is $\mc{Z}$-stable and $P$ is nuclear, then $P$ is
$\mc{Z}$-stable.
\item If $A$ is stably finite and exact and has strict comparison, then $P$
has strict comparison.
\item If $A$ is stably finite and has tracial $m$-comparison, then $P$ has
tracial $m$-comparison.
\item If $A$ is exact and tracially $m$-almost divisible, then $P$ is
tracially $m$-almost divisible.
\item If $A$ and $P$ are nuclear and
$\operatorname{Trdim}_{\mathrm{nuc}}(A)\leq m$, then
$\operatorname{Trdim}_{\mathrm{nuc}}(P)\leq m$.
\end{enumerate}
Here $m\in\mathbb{N}$ in (4) and (5), and $m\in\mathbb{N}\cup\{0\}$ in
(6).
\end{thm}

\begin{proof}
Proposition~\ref{P:inclusion-sequentiallysplit} shows that
$\iota:P\hookrightarrow A$ is tracially sequentially-split by order zero;
moreover, $P$ is simple by Proposition~\ref{P:Rokhlinsimple}. Assertions (1)
and (2) now follow from \cite[Theorem 3.6]{LeeOsaka4}, and (3) follows from
\cite[Theorem 3.10]{LeeOsaka4}. Assertions (4), (5), and (6) follow,
respectively, from \cite[Theorems 3.3, 3.6, and 3.10]{LeePermanenceII}.
\end{proof}

\begin{cor}\label{C:framework-stable-rank}
Under the assumptions of Theorem~\ref{T:framework-permanence}, if $A$ is
stably finite and $\mc{Z}$-stable and both $A$ and $P$ are nuclear, then
$\tsr(P)=1$.
\end{cor}

\begin{proof}
This is \cite[Corollary 3.7]{LeeOsaka4}, applied to the inclusion map in
Proposition~\ref{P:inclusion-sequentiallysplit}.
\end{proof}

The same framework also packages the consequences of the duality theorem.

\begin{cor}\label{C:duality}
Let $P\subset A$, $E:A\to P$, and
$B=C\sp*\langle A,e_P\rangle$ satisfy the assumptions of
Theorem~\ref{T:tracialduality2}, and suppose that $B$ is simple and finite. If
$E$ is weakly tracially approximately representable, then the inclusion
$A\hookrightarrow B$ is tracially sequentially-split by order zero.
Consequently, each of the conclusions in
Theorem~\ref{T:framework-permanence}, with $(P,A)$ replaced by $(A,B)$ and
with the corresponding hypotheses, passes from $B$ to $A$.
\end{cor}

\begin{proof}
By Theorem~\ref{T:tracialduality2}, the dual conditional expectation
$\widehat{E}:B\to A$ has the weak tracial Rokhlin property. Applying
Proposition~\ref{P:inclusion-sequentiallysplit} to $A\subset B$ and
$\widehat{E}$ gives the first assertion, and
Theorem~\ref{T:framework-permanence} gives the rest.
\end{proof}

We finish by recording the finite-group case in the same language.

\begin{cor}\label{C:action-framework}
Let $G$ be a finite group and let $\alpha:G\curvearrowright A$ be an action
with the weak tracial Rokhlin property on a simple, separable, unital, finite,
infinite-dimensional $C\sp*$-algebra. Then the canonical inclusion
$A^{\alpha}\hookrightarrow A$ is tracially sequentially-split by order zero,
and the conclusions of Theorem~\ref{T:framework-permanence} pass from $A$ to
$A^{\alpha}$ under their respective hypotheses.

Let $\sigma:G\curvearrowright C(G)$ denote the left translation action,
defined by
\[
   (\sigma_g(f))(h)=f(g^{-1}h),\qquad f\in C(G),\quad g,h\in G.
\]
Then the canonical equivariant homomorphism
$1_{C(G)}\otimes\id_A:(A,\alpha)\to
(C(G)\otimes A,\sigma\otimes\alpha)$ induces a homomorphism
\[
(1_{C(G)}\otimes\id_A)\rtimes G:
A\rtimes_{\alpha}G\longrightarrow
(C(G)\otimes A)\rtimes_{\sigma\otimes\alpha}G
\cong M_{|G|}(A)
\]
that is tracially sequentially-split by order zero. Hence the same permanence
theorems apply from $M_{|G|}(A)$ to $A\rtimes_{\alpha}G$ whenever their
respective hypotheses are satisfied.
\end{cor}

\begin{proof}
Theorem~\ref{T:action-inclusion} and
Proposition~\ref{P:inclusion-sequentiallysplit} give the assertion for the
fixed-point algebra. The crossed-product homomorphism is tracially
sequentially-split by order zero by \cite[Theorem 4.8]{LeeOsaka4}; applying
\cite[Theorems 3.6 and 3.10]{LeeOsaka4} and
\cite[Theorems 3.3, 3.6, and 3.10]{LeePermanenceII} gives the final assertion.

For completeness, we explain the matrix-algebra identification, which does
not require $G$ to be abelian. Under the identification
$C(G)\otimes A\cong C(G,A)$, define
\[
   \Theta:C(G,A)\longrightarrow C(G,A),\qquad
   (\Theta(F))(h)=\alpha_{h^{-1}}(F(h)).
\]
Then
\[
   \Theta\circ(\sigma_g\otimes\alpha_g)
   =(\sigma_g\otimes\id_A)\circ\Theta
   \qquad (g\in G).
\]
Consequently,
\[
\begin{aligned}
(C(G)\otimes A)\rtimes_{\sigma\otimes\alpha}G
&\cong (C(G)\otimes A)\rtimes_{\sigma\otimes\id_A}G\\
&\cong (C(G)\rtimes_\sigma G)\otimes A\\
&\cong \mathcal{K}(\ell^2(G))\otimes A
 \cong M_{|G|}(A).
\end{aligned}
\]
Thus the assertion uses the translation action on $C(G)$ rather than
Pontryagin duality and holds for every finite group.
\end{proof}

\section{Acknowledgements}
This research was initiated during the first author's visit to the University
of Oregon in 2017. He thanks the host, N.~Christopher Phillips, for pointing
out the orthogonalization argument used in
Lemma~\ref{L:ultrapowerorthogonalelements}. Part of this work was also
completed during the first author's visit to Ritsumeikan University. He
gratefully acknowledges the support of the Department of Mathematical
Sciences at Ritsumeikan University.



\begin{thebibliography}{99}

\bibitem{AsadiVasfi} M. A. Asadi-Vasfi, \emph{Weakly tracially approximately representable actions}, J. Operator Theory \textbf{91} (2024), no. 1, 3--25.

\bibitem{BS} S. Barlak and G. Szab\'{o}, \emph{Sequentially split $*$-homomorphisms between $C\sp*$-algebras}, Internat. J. Math. \textbf{27} (2016), no. 12, 1650105, 48 pp.

\bibitem{Cu} J. Cuntz, \emph{Dimension functions on simple $C\sp*$-algebras}, Math. Ann. {\bf 233} (1978), 145--153.

\bibitem{ET} G. A. Elliott and A. Toms, \emph{Regularity properties in the classification program for separable amenable $C\sp*$-algebras}, Bull. Amer. Math. Soc. \textbf{45} (2008), no. 2, 229--245.

\bibitem{EvansKawahigashi} D. E. Evans and Y. Kawahigashi,
\emph{Quantum Symmetries on Operator Algebras}, Oxford Mathematical
Monographs, The Clarendon Press, Oxford University Press, New York, 1998.

\bibitem{Ha} U. Haagerup, \emph{Quasitraces on exact $C\sp*$-algebras are traces}, C. R. Math. Acad. Sci. Soc. R. Can. \textbf{36} (2014), nos. 2--3, 67--92.

\bibitem{G} E. Gardella, \emph{Crossed products by compact group actions with the Rokhlin property}, J. Noncommut. Geom. \textbf{11} (2017), no. 4, 1593--1626.

\bibitem{GHS} E. Gardella, I. Hirshberg, and L. Santiago, \emph{Rokhlin dimension: duality, tracial properties, and crossed products}, Ergodic Theory Dynam. Systems \textbf{41} (2021), 408--460.

\bibitem{GHV} E. Gardella, I. Hirshberg, and A. Vaccaro, \emph{Strongly outer actions of amenable groups on $\mathcal{Z}$-stable nuclear $C\sp*$-algebras}, J. Math. Pures Appl. \textbf{162} (2022), 76--123.



\bibitem{HW:Rokhlin} I. Hirshberg and W. Winter,  \emph{Rokhlin actions and self-absorbing $C\sp*$-algebras}, Pacific J. Math. \textbf{233} (2007), no. 1,  125--143.

\bibitem{Izumi:finite} M. Izumi, \emph{ Finite group actions on $C\sp*$-algebras with the Rohlin property I}, Duke Math. J. \textbf{122} (2004), no. 2,  233--280.

\bibitem{Izumi:Rokhlin} M. Izumi, \emph{Finite group actions on $C\sp*$-algebras with the Rohlin property II}, Adv. Math. \textbf{184} (2004), 119--160.


\bibitem{Izumi:inclusion} M. Izumi, \emph{Inclusions of simple $C\sp*$-algebras}, J. Reine Angew. Math. \textbf{547} (2002), 97--138.

\bibitem{JiangSu} X. Jiang and H. Su, \emph{On a simple unital projectionless $C\sp*$-algebra}, Amer. J. Math. \textbf{121} (1999), no. 2, 359--413, doi:10.1353/ajm.1999.0012.



\bibitem{JP:saturated} J. Jeong and  K. Park,  \emph{Saturated actions by finite-dimensional Hopf *-algebras on $C\sp*$-algebras}, Internat. J. Math. \textbf{19}(2008), no.2, 125--144. 

\bibitem{JOPT:cancellation} J. A. Jeong, H. Osaka, N. C. Phillips, and T. Teruya,
\emph{Cancellation for inclusions of $C\sp*$-algebras of finite depth},
Indiana Univ. Math. J. \textbf{58} (2009), no. 4, 1537--1564.



\bibitem{KR:positive} E. Kirchberg and M. R{\o}rdam, \emph{Non-simple purely infinite $C\sp*$-algebras}, Amer. J. Math. \textbf{122}(2000), no.3, 637--666. 

\bibitem{KW:lifting} E. Kirchberg and  W. Winter,  \emph {Covering dimension and quasidiagonality}, Internat. J. Math. {\bf 15} (2004), no. 1,  63--85.

\bibitem{Kishi1}A. Kishimoto,  \emph{A Rohlin property for one-parameter automorphism groups}, Comm. Math. Phys. \textbf{179}(1996), no. 3,  599--622.

\bibitem{Kishi2}A. Kishimoto \emph{Outer Automorphisms and Reduced Crossed Products of simple $C\sp*$-algebras},  Comm. Math. Phys. \textbf{81}(1981), 429--435. 



\bibitem{Lin:classification}H. Lin,  \emph{Classification of simple $C\sp*$-algebras of tracial topological rank zero},  Duke Math. J. {\bf125}(2004), 91--119.

\bibitem{LeeOsaka} H. Lee and H. Osaka, \emph{Tracially sequentially-split $*$-homomorphisms between $C\sp*$-algebras}, Ann. Funct. Anal. {\bf 12} (2021), Article 37.

\bibitem{LO:2019} H. Lee and H. Osaka, \emph{On dualities of actions and inclusions}, J. Funct. Anal. \textbf{276} (2019), 602--635.

\bibitem{LOT:Rokhlin-index} H. Lee, H. Osaka, and T. Teruya,
\emph{The Rokhlin inclusions with integer and non-integer index},
arXiv:2405.03964v2 [math.OA] (2025).

\bibitem{LeeOsaka4} H. Lee and H. Osaka, \emph{On permanence of regularity properties}, J. Topol. Anal. \textbf{17} (2025), no. 1, 35--54.

\bibitem{LeePermanenceII} H. Lee, \emph{On permanence of regularity properties II}, arXiv:2603.07491 (2026).

\bibitem{Lee} H. Lee, \emph{On dualities of actions}, J. Math. Anal. Appl. \textbf{537} (2024), no. 2, Article 128268.

\bibitem{MS12} H. Matui and Y. Sato,  \emph{$\mathcal{Z}$-stability of crossed products by strongly outer actions}, Comm. Math. Phys.  {\bf 314}  (2012),  no. 1, 193--228.

\bibitem{Nw16} N. Nawata,  \emph{Finite group actions on certain stably projectionless $C\sp*$-algebras with the Rohlin property}, Trans. Amer. Math. Soc. \textbf{368} (2016), 471--493.


\bibitem{Osaka:(SP)-property} H. Osaka, \emph{(SP)-property for a pair of $C\sp*$-algebras}, J. Operator Theory {\bf 46} (2001), 159--171.

\bibitem{OT:tsr} H. Osaka and T. Teruya,
\emph{Topological stable rank of inclusions of unital $C\sp*$-algebras},
Internat. J. Math. \textbf{17} (2006), no. 1, 19--34.

\bibitem{OKT:Rokhlin} H. Osaka, K. Kodaka, and T. Teruya, \emph{Rokhlin property for inclusion of $C\sp*$-algebras with a finite Watatani index}, in \emph{Operator Structures and Dynamical Systems}, Contemp. Math. \textbf{503}, Amer. Math. Soc., Providence, RI, 2009, 177--195.

\bibitem{OT1} H. Osaka and T. Teruya,  \emph{The Jiang-Su absorption for inclusions of $C\sp*$-algebras},  Canad. J. Math.  {\bf 70} (2018), no.2, 400--425.

\bibitem{OT2} H. Osaka and T. Teruya, \emph{Erratum: The Jiang--Su absorption for inclusions of unital $C\sp*$-algebras}, Canad. J. Math. {\bf 73} (2021), no. 1, 293--295, doi:10.4153/S0008414X20000310.

\bibitem{OP06} H. Osaka and N. C. Phillips, \emph{Stable and real rank for crossed products by automorphisms with the tracial Rokhlin property}, Ergodic Theory Dynam. Systems  {\bf26}  (2006),  no. 5, 1579--1621.

\bibitem{Phillips:tracial} N. C. Phillips, \emph{The tracial Rokhlin property for  actions of finite groups on $C\sp*$-algebras}, Amer. J. Math. \textbf{133}(2011), no. 3, 581--636.

\bibitem{Phillips:large} N. C. Phillips, \emph{Large subalgebras}, arXiv:1408.5546v1.

\bibitem{Ro:UHF1} M. R\o rdam, \emph{On the structure of simple $C\sp*$-algebras tensored with a UHF-algebra}, J. Funct. Anal. \textbf{100} (1991), 1--17.

\bibitem{Ro:UHF2} M. R\o rdam, \emph{On the structure of simple $C\sp*$-algebras tensored with a UHF-algebra II}, J. Funct. Anal. \textbf{107} (1992), 255--269.

\bibitem{SP:saturated} W. Szyma\'nski and C. Peligrad, \emph{Saturated actions of finite dimensional Hopf $*$-algebras on $C\sp*$-algebras}, Math. Scand. \textbf{75}(1994), 217--239.

\bibitem{Szlachanyi} K. Szlach\'anyi, \emph{Weak Hopf algebra symmetries of $C\sp*$-algebra inclusions}, arXiv:math/0101005 (2001).

\bibitem{Ro:Classification} M. R\o rdam and E. St\o rmer, \emph{Classification of Nuclear $C\sp*$-Algebras: Entropy in Operator Algebras}, Encyclopaedia of Mathematical Sciences, vol. \textbf{126}, Springer-Verlag, Berlin, 2002.

\bibitem{San15}  L. Santiago, \emph{Crossed products by actions of finite groups with the Rokhlin property},
Internat. J. Math. 26 (2015), no. 7, 1550045.

\bibitem{TW:ssa} A. Toms and W. Winter, \emph{Strongly self-absorbing $C\sp*$-algebras}, Trans. Amer. Math. Soc. \textbf{359} (2007), no. 8, 3999--4029.

\bibitem{W:index} Y. Watatani, \emph{Index for $C\sp*$-Subalgebras}, Mem. Amer. Math. Soc. \textbf{83} (1990), no. 424.

\bibitem{WZ:orderzero} W. Winter and J. Zacharias, \emph{Completely positive maps of order zero}, M\"unster J. Math. {\bf 2} (2009), 311--324.

\enlargethispage{8\baselineskip}
\end{thebibliography}
\end{document}